\documentclass[12pt]{amsart}
\usepackage{amsbsy}
\usepackage{graphicx,epsfig,subfigure,psfrag}
\begin{document}

\newtheorem{theorem}{Theorem}
\newtheorem{proposition}{Proposition}
\newtheorem{lemma}{Lemma}
\newtheorem{corollary}{Corollary}
\newtheorem{definition}{Definition}
\newtheorem{remark}{Remark}
\newcommand{\beq}{\begin{equation}}
\newcommand{\eeq}{\end{equation}}
\numberwithin{equation}{section} \numberwithin{theorem}{section}
\numberwithin{proposition}{section} \numberwithin{lemma}{section}
\numberwithin{corollary}{section}
\numberwithin{definition}{section} \numberwithin{remark}{section}
\newcommand{\ren}{\mathbb{R}^N}
\newcommand{\re}{\mathbb{R}}
\newcommand{\n}{\nabla}
\newcommand{\iy}{\infty}
\newcommand{\pa}{\partial}
\newcommand{\fp}{\noindent}
\newcommand{\ms}{\medskip\vskip-.1cm}
\newcommand{\mpb}{\medskip}
\newcommand{\BB}{{\bf B}}
\newcommand{\Am}{{\bf A}_{2m}}
\renewcommand{\a}{\alpha}
\renewcommand{\b}{\beta}
\newcommand{\g}{\gamma}
\newcommand{\G}{\Gamma}
\renewcommand{\d}{\delta}
\newcommand{\D}{\Delta}
\newcommand{\e}{\varepsilon}
\newcommand{\var}{\varphi}
\renewcommand{\l}{\lambda}
\renewcommand{\o}{\omega}
\renewcommand{\O}{\Omega}
\newcommand{\s}{\sigma}
\renewcommand{\t}{\tau}
\renewcommand{\th}{\theta}
\newcommand{\z}{\zeta}
\newcommand{\wx}{\widetilde x}
\newcommand{\wt}{\widetilde t}
\newcommand{\noi}{\noindent}
\newcommand{\inA}{\quad \mbox{in} \quad \ren \times \re_+}
\newcommand{\inB}{\quad \mbox{in} \quad}
\newcommand{\inC}{\quad \mbox{in} \quad \re \times \re_+}
\newcommand{\inD}{\quad \mbox{in} \quad \re}
\newcommand{\forA}{\quad \mbox{for} \quad}
\newcommand{\whereA}{,\quad \mbox{where} \quad}
\newcommand{\asA}{\quad \mbox{as} \quad}
\newcommand{\andA}{\quad \mbox{and} \quad}
\newcommand{\ssk}{\smallskip}
\newcommand{\LongA}{\quad \Longrightarrow \quad}
\def\com#1{\fbox{\parbox{6in}{\texttt{#1}}}}

\title
{\bf
Nonlinear dispersion equations:\\ smooth deformations,
compactons,\\ and extensions to higher orders}



\author {
Victor A.~Galaktionov}

\address{Department of Mathematical Sciences, University of Bath,
 Bath BA2 7AY, UK}
\email{vag@maths.bath.ac.uk}


\thanks{Research partially supported by the INTAS network
CERN-INTAS00-0136}

  \keywords{Odd-order quasilinear PDE, shock and rarefaction waves,
  entropy solutions,  self-similar patterns.
  {\bf To appear in:}  Comp. Math. Math. Phys.
  }
 \date{\today}

\begin{abstract}

 The third-order nonlinear
dispersion PDE, as the key model,
 \beq
 \label{01}
 u_t=(uu_x)_{xx} 
 \quad \mbox{in}
 \quad \re \times \re_+,
 \eeq
 is studied.
 Two  Riemann's problems
  for (\ref{01})
 with initial data
  $
  S_\mp(x) = \mp{\rm sign}\, x,
  $
   create
    the shock ($u(x,t) \equiv S_-(x)$) and smooth rarefaction (for data $S_+$) waves, \cite{3NDEI}.

The  concept of ``$\d$-entropy" solutions (a``$\d$-entropy test")
and others are developed for distinguishing
 shock and rarefaction waves
 by using stable smooth $\d$-deformations
of discontinuous solutions.
These are analogous to entropy solutions for
   scalar conservation laws such as
 $
 u_t + uu_x=0,
  $
 developed by Oleinik and Kruzhkov
 (in $\ren$)
in the 1950-60s.
The Rosenau--Hyman $K(2,2)$ (compacton) equation
 $$
 u_t=(u u_x)_{xx} + 4 uu_x,
 $$
 which has a special importance for applications, is studied.
 Compactons as compactly supported travelling wave solutions
 are shown to pass the $\d$-entropy test.
  Shock and rarefaction waves
  are discussed  for   other NDEs
  such as
 $$
 u_t=(u^2 u_x)_{xx}, \,\,\,  u_{tt}=(u u_x)_{xx}, \,\,\,u_{tt}=
 uu_x, \,\,\,
  u_{ttt}=(u u_x)_{xx}, \,\,\, u_t=(uu_x)_{xxxxxx},
\,\,\mbox{etc}.
 $$

\end{abstract}

\maketitle


\begin{center}
{\em Dedicated to  the memory of Professors O.A.~Oleinik and
S.N.~Kruzhkov}
\end{center}




\section {Introduction: nonlinear dispersion PDEs and main results}
 \label{Sect1}

 \subsection{NDEs: nonlinear dispersion equations in application and general PDE theory}

The present paper continues the study began in \cite{3NDEI}
 of odd-order {\em nonlinear
dispersion} (or {\em dispersive}) PDEs (NDEs).
The canonical  model 
is the 
{\em third-order quadratic NDE} (the NDE--$3$)
 \beq
 \label{1}
  \mbox{$
 u_t={\bf A}(u) \equiv (uu_x)_{xx} = u u_{xxx}+ 3 u_x u_{xx}
 \quad \mbox{in}
 \quad \re \times (0,T), \,\,\, T >0.
  $}
 \eeq
Posing for (\ref{1}) the Cauchy problem
  includes locally
integrable initial data
 \beq
 \label{2}
 u(x,0)=u_0(x) \quad \mbox{in} \quad \re.
  \eeq
  Frequently, we assume that $u_0$ is bounded and compactly supported.
We will also deal with
the initial-boundary values problem in $(-L,L)\times \re_+$ with
Dirichlet boundary conditions. Main applications concerning NDEs
can be found in \cite{GalNDE5, 3NDEI}; see also
\cite[Ch.~4]{GSVR}, so that we do not discuss these issues in
detail. However, we need to stress the attention of the  Reader to
the compacton phenomena, which were not properly treated in the
mathematical literature.

\ssk

\noi\underline{\em Compact patterns and NDEs}.
 These are known for the  {\em Rosenau--Hyman} (RH)
{\em equation}
  \beq
  \label{Comp.4}
  \mbox{$
  u_t =  (u^2)_{xxx} + (u^2)_x,
  $}
  \eeq
   which is the $K(2,2)$ equation from
    the general $K(m,n)$ family of
   the NDEs:
  \beq
  \label{Comp.5}
   u_t =  (u^n)_{xxx} +  (u^m)_x \quad (u \ge 0).
   \eeq
 References on physical applications  of such NDEs
 are available in
  \cite[\S~1]{3NDEI} and in \cite[\S~4.2]{GSVR}.
  We will check entropy properties of  compactons
 for various NDEs of this type.



Further applied compacton-like models are discussed in
\cite{3NDEI}. A standard definition of weak solutions for
(\ref{1}) is also presented there, so that we are in a position to
explain our main targets concerning entropy-like theory of shocks.

\subsection{Plan of the paper: entropy theory (a test) via smooth
deformations and compactons}
 \label{Fok1}

As in \cite{3NDEI}, we begin with discussion of  some auxiliary
properties of the NDE--3.

\ssk

\noi\underline{\em Smoothing for the NDE--3}. Firstly, we recall
that the smoothing phenomena and results for sufficiently regular
solutions of linear and nonlinear third-order PDEs are well know
from the 1980-90s. For instance, infinite $C^\infty$-smoothing
results were proved in \cite{Cr90} for a general linear equation
of the form
 \beq
 \label{lin11}
 u_t+a(x,t) u_{xxx}=0 \quad (a(x,t) \ge  c>0),
 \eeq
 and in \cite{Cr92} for the corresponding fully nonlinear PDE
  \beq
  \label{Lin12}
  u_t+f(u_{xxx},u_{xx},u_x,u,x,t)=0 \quad \big(\,f_{u_{xxx}} \ge
  c>0\, \big).
   \eeq
   Namely, for a class of such equations, it is shown that, for data
   with minimal regularity and sufficient decay at infinity, there
   exists a unique
   solution $u(x,t) \in C^\infty$ for arbitrarily small $t>0$.
 Similar smoothing local in time results for unique solutions
 are available for
  \beq
  \label{lin13}
  u_t + f(D^3u,D^2u,Du,u,x,y,t)=0 \quad \mbox{in} \quad \re^2 \times \re_+;
   \eeq
   see \cite{Lev01} and further  references therein. Concerning
unique continuation and continuous dependence properties, see
\cite{Daw07} and references therein, and \cite{Tao00} for various
estimates.

\ssk

\noi\underline{\em The NDE: a conservation law in $H^{-1}$}.
Writing (\ref{1}) as (see details in \cite[\S~1.4]{3NDEI})
 \beq
 \label{10}
 (-D_x^2)^{-1}u_t + uu_x=0,
  \eeq
 yields the first {\em a priori} uniform bound
 for data $u_0 \in H^{-1}(\re)$. Namely,  multiplying (\ref{10}) by $u$ in
$L^2$ gives the conservation law
 \beq
 \label{111}
  \mbox{$
 \frac 12 \, \frac{\mathrm d}{{\mathrm d}t}\, \|u(t)\|_{H^{-1}}^2
 =0
 \quad \Longrightarrow \quad \|u(t)\|_{H^{-1}}^2= c_0= \|u_0\|_{H^{-1}}^2
 \quad \mbox{for all} \,\,\, t>0.
  $}
  \eeq


\ssk

\noi\underline{\em Main results}.
 In the present
paper, we propose some concepts for developing
adequate mathematics of NDEs with shocks, which will be concluded
in Section \ref{Sect9} by revealing connections with other classes
of nonlinear degenerate PDEs. It turns out that some NDE concepts
has definite reliable common roots and can be put into the
framework of much better developed theory of quasilinear parabolic
equations. We restrict our attention to a key demand, how to
distinguish the shock and rarefaction waves, and this is done by
developing the so-called "$\d$-entropy test" on solutions via
smooth deformations. General uniqueness-entropy theory for NDEs
such as (\ref{1}) and others is shown to be illusive
\cite{Gal3NDENew}.

 Concerning the simple canonical  model (\ref{1}),
 we  do
 the following:


 $\bullet$ {\bf (i)} Reviewing local existence and uniqueness theory for
 the NDE (\ref{1}) and, on its basis, developing an  $\d$-entropy test
 for distinguishing shock and rarefaction waves.


 For the RH equation such as (\ref{Comp.4}), we prove that:


$\bullet$  {\bf (ii)} Rosenau's compacton solutions are both
$\d$-entropy and G-admissible.


 Some of related  questions and results were previously discussed in
 a more applied and formal fashion in \cite[\S~7]{GalEng} and
 \cite[Ch.~4]{GSVR}.

\subsection{On extensions and other nonlinear dispersion models}

The developed concepts  cover a wide range of various NDEs. First
of all, we should mention that the fact that (\ref{1}) is
degenerate at $u=0$ and hence admits compactly supported solutions
(which is an interesting pleasant feature) makes the analysis of
$\d$-entropy solutions and shocks much harder.
 However, shock waves exist for other non-degenerated NDEs with
 analytic coefficients. For instance, we study entropy
 shocks for the NDE with infinite propagation,
  \beq
 \label{non1}
 u_t=((1+u^2)u_x)_{xx}.
  \eeq
All
our further NDEs admit analogous non-degenerate versions
admitting shock and rarefaction waves, but no finite propagation
and interfaces in the Cauchy problem.

Another  related to (\ref{1}) model to be  discussed is
the {\em cubic fully divergent NDE}
 \beq
 \label{HD1N}
  \mbox{$
  u_t=(u^2 u_x)_{xx} \equiv \frac 13 \, (u^3)_{xxx} \quad \big(\mbox{the conservation law analogy is}
  \,\,\, u_t + u^2 u_x=0\big).
   $}
   \eeq
   We study (\ref{HD1N}) instead of less physically motivated
``quadratic" model
 $
 u_t=(|u|u_x)_{xx}
 $
 that exhibits similar properties of shocks and rarefaction waves.

The results on $\d$-entropy solutions and similarity patterns can
be extended (Section \ref{Sect9}) to truly quadratic {\em
non-fully divergent} NDEs such as
 \beq
 \label{3,1}
 u_t=(u u_{xx})_x \equiv u u_{xxx} + u_x u_{xx},
  \eeq
  which we call the NDE--(2,1), where 2 and 1 stand for the number
  of the internal and external derivatives in this differential
  form. Notice that a standard concept of weak solutions hardly
  applies to (\ref{3,1}), so that the shock $S_-(x)$ is not a weak
  solution.
  In order to underline once more the fact that being
  weak is not a necessary demand, we consider a formal {\em
  fully nonlinear NDE}
   \beq
   \label{FN1}
   |u_t|^\g u_t=(uu_x)_{xx}, \quad \mbox{where} \quad \g >-1.
    \eeq
For $\g=0$, this gives the original equation (\ref{1}). Obviously,
for $\g \not = 0$, (\ref{FN1}) does not admit any weak
formulation. Nevertheless, we show that (\ref{FN1}) admits blow-up
formation of shocks of $S_-$-type.

  In Section \ref{Sect2t}, we discuss the
 shock
  formation mechanism for higher-order in time NDEs,
   \beq
   \label{222}
   u_{tt}=(uu_x)_{xx} \quad \mbox{and}
   \quad  u_{ttt}=(uu_x)_{xx}.
    \eeq

 Several principal features
remain the same for higher-order NDEs such as the quadratic {\em
fifth-order NDE} (NDE--5)
 \beq
 \label{N5}
u_t=-(u u_x)_{xxxx}
 \quad
\mbox{or, in general,} \quad u_t=(-1)^{m+1}D_x^{2m}(u u_x), \,\,\,
m \ge 1;
 \eeq
 see Section \ref{Sect55}.
 These are conservation laws in $H^{-2}$, or $H^{-m}$.
The mathematics of particular similarity solutions with shocks is
developed in similar lines but technically becomes more involved,
so we have to  catch the similarity profiles numerically.

We also claim that some concepts such as smooth $\d$-deformation
 and others,
 developed for models in 1D can be
adapted to the $N$-dimensional NDEs. In particular, the basic NDE
(\ref{1}) in $\ren$ takes the form
 \beq
 \label{ndeN}
  \mbox{$
 u_t = \D (u\frac{\partial u}{\partial x_1}) \equiv \frac 12\,
 \frac{\partial}{\partial x_1}\, \D u^2 \quad \mbox{in} \quad
 \ren\times (0,T).
 $}
 \eeq


\section{Conservation laws: smooth $\d$-deformations define
entropy solutions}
 \label{Sect40}

From now on, being sufficiently informed about formation of
crucial shock and other singularities in the NDEs, we will start
to investigate the general questions on existence and uniqueness
of entropy weak solutions of  (\ref{1}).
As usual, we begin our discussion by stressing attention to key
analogies with classic theory of  first-order conservation laws
such as Euler's equation from gas dynamics
 \beq
 \label{3}
  u_t + uu_x=0
\quad \mbox{in} \quad \re\times \re_+.
 \eeq
 Entropy theory for such first-order PDEs  was created by Oleinik \cite{Ol1, Ol59} and Kruzhkov
\cite{Kru2} ($x \in \ren$) in the 1950--60s; see details on the
history, main results, and modern developments in the well-known
monographs \cite{Bres, Daf, Sm}.
 Thus,  we now apply smooth  $\d$-deformation
concepts to these simpler PDEs considered now  in $Q_1= \re \times
(0,1)$.

\subsection{Preliminaries: entropy inequalities and  solutions for conservation laws}
\label{S6.1}


It is known from  the 1950's  that the Cauchy problem for general
scalar conservation laws admits a unique entropy solution. We
refer to first complete results by Oleinik (obtained in 1954-56),
who introduced entropy conditions  in 1D and proved existence and
uniqueness results (see survey \cite{Ol1}), and by Kruzhkov (1970)
\cite{Kru2}, who developed general non-local theory of entropy
solutions in $\ren$. In the general case, one of Oleinik's local
entropy condition has the form \cite[p.~106]{Ol1}
 \beq
 \label{1.2L}
  \mbox{$
  \frac{u(x_1,t) - u(x_2,t)}{x_1 -x_2} \le K(x_1,x_2,t) \quad
  \mbox{for all} \,\,\, x_1, x_2 \in \re, \,\, t \in [0,1],
   $}
 \eeq
 where $K$ is a continuous function for $t \in [0,1]$.
  Oleinik's local condition E (Entropy) introduced in
 \cite{Ol59},
 for the model equation (\ref{3})
 corresponds to the well-known principle of non-increasing  entropy
 from gas dynamics,
 \beq
 \label{Ent2}
 u(x^+,t) \le u(x^-,t) \quad \mbox{in} \,\,\, Q_1= \re \times (0,1],
 \eeq
 with strict inequality on lines of discontinuity, \cite[p.~101]{Ol1}.

{\em Kruzhkov's entropy
 condition} \cite{Kru2}
on  solutions $u \in L^\infty(Q_1)$ of (\ref{3})  takes the form
of
 the non-local inequality
 \beq
 \label{KrIn}
 \mbox{$
 |u-k|_t + \frac 12 \, \big[{\rm sign}(u-k)  (u^2-k^2)\big]_x \le 0 \quad \mbox{in} \,\,\,
 {\mathcal D}'(Q_1) \quad \mbox{for any} \,\,\, k \in \re.
  $}
  \eeq
  This inequality  is understood in the sense of distributions
 meaning that the sign $\le$ is preserved after
  multiplying the inequality by any smooth compactly supported cut-off
  function $\varphi \in C_0^\infty(Q_1)$,
    $\varphi \ge 0$,
  and integrating by parts. See  clear presentation of these ideas
 in
 Taylor \cite[p.~401]{Tay}.
Oleinik's and Kruzhkov's approaches are known to coincide in the
1D geometry. Both entropy conditions generate a semigroup of
contractions in $L^1$, so that if $u$ and $v$ are two solutions of
(\ref{3}), then
 \beq
 \label{sem1}
 \mbox{$
  \frac{\mathrm d}{{\mathrm d}t} \, \|u(t)- v(t) \|_{L^1} \le 0.
   $}
    \eeq

It is  key that the unique entropy solution is constructed by the
parabolic $\e$-approximation
 \beq
 \label{cle1}
 u_\e: \quad u_t+ u u_x = \e u_{xx} \quad (\e>0).
  \eeq
 Multiplying (\ref{cle1})
 by  any smooth
monotone increasing function $E(u)$ (an approximation of ${\rm
sign} \,(u-k)$ for any  $k \in \re$) yields on integration by
parts the correct sign:
 \beq
 \label{cc1}
  \mbox{$
\iint \e u_{xx} E(u)=- \e \iint E'(u)(u_x)^2 \le 0.
 $}
 \eeq
Hence, as $\e \to 0$, this gives the necessary sign as in
(\ref{KrIn}).


The obvious advantage of the conservation law (\ref{3})
 is that, for smooth initial data (\ref{2}), the unique local continuous solution
is obtained by method of characteristics and is given by the
corresponding algebraic equation
 \beq
 \label{p1}
  \mbox{$
   {\mathrm d}t= \frac{{\mathrm d}x}u \quad \Longrightarrow \quad
u(x,t) = u_0(x- u(x,t)t) \quad \mbox{for all} \quad t \in [0,\D
t),
 $}
 \eeq
 where $\D t \le 1$ is the first moment of time when  a
 shock of the type $S_-(x)$ (this type is guaranteed by (\ref{Ent2}))
 occurs at some point or many points.


Thus,
 for $t \ge \D t$, it is necessary to apply the entropy
 inequalities to select good (entropy) solutions.
Using this, and bearing in mind that entropy solutions are
continuous relative initial data (in $L^1$, say), we propose the
following construction which is fully based on algebraic relations
(\ref{p1}):

\subsection{Conservation laws: $\d$-stable ${\mathbf=}$ entropy solutions}
 \label{S6.2}

It is the obvious well-known and, nevertheless, crucial
observation that, by the characteristic mechanism  (\ref{p1}),
 \beq
 \label{p23}
 \mbox{
 non-entropy shocks of the shape $S_+$ cannot appear
 evolutionary.}
 \eeq
Indeed, 
 differentiating (\ref{p1}) in $x$
yields
 \beq
 \label{p33}
  \begin{matrix}
  \mbox{$
 u_x(x,t)= \frac {u_0'(x-u(x,t)  t)}{1+u_0'(x-u(x,t) t)t}, \quad \mbox{so that}
  $}\qquad\qquad \ssk\\
u_0' \ge 0 \quad \Longrightarrow \quad \mbox{no blow-up of $u_x$
(``gradient catastrophe") occurs.}\qquad\qquad
 \end{matrix}
 \eeq

Recalling the necessary evolution property in (\ref{p33}), given a
small $\d>0$ and a bounded (say, for simplicity, in $L^1$ and in
$L^\infty$) solution $u(x,t)$ of the Cauchy problem (\ref{3}),
(\ref{2}), we construct its  {\em $\d$-deformation} given
explicitly by the characteristic method (\ref{p1}) as follows:


(i) we perform a smooth $\d$-deformation of initial data $u_0 \in
L^1\cap L^\infty$ by introducing a suitable $C^1$ function
$u_{0\d}(x)$ such that
 \beq
 \label{p2}
 \mbox{$
\int |u_0-u_{0\d}| < \d.
 $}
 \eeq
 By $u_{1\d}(x,t)$ we denote the unique local solution of the Cauchy
 problem with data $u_{0\d}$, so that by (\ref{p1}), continuous
 function
 $u_{1\d}(x,t)$ is defined algebraically on the maximal interval $t \in
 [t_0,t_1(\d))$, where we denote $t_0=0$ and  $t_1(\d)=\D_{1\d}$.
 It is important that, here and later on, smooth deformations are performed in a
 small neighbourhood of possible discontinuities {\sc only}
 leaving the rest of smooth profiles untouchable, so that these evolve along the characteristics,
 as usual.

 Actually, this emphasizes the obvious fact that the shocks (on a
 set of zero measure) occur as
 a result of nonlinear interaction of the areas with continuous
 solutions, which hence cannot be connected without
 discontinuities.


 (ii) Since at $t= \D_{1\d}$ a shock of type $S_-$
  (or possibly infinitely many shocks) is supposed to
 occur, since otherwise we continue the algebraic procedure, we
 perform another suitable $\d$-deformation of the ``data"
 $u_{1\d}(x,\D_{1\d})$ to get a unique continuous solution $u_{2\d}(x,t)$
 on the maximal interval $t \in [t_1(\d),t_2(\d))$, with
 $t_2(\d)=\D_{1\d}+\D_{2\d}$, etc.


$\dots$


 (k) 
 With  suitable
choices of each $\d$-deformations of ``data" at the moments
$t=t_j(\d)$, when $u_{j\d}(x,t)$ has a shock for $j=1,2,...$,
there exists a $t_{k}(\d)
> 1$ for some finite $k=k(\d)$, where $k(\d) \to +\infty$ as $\d
\to 0$. It is easy to see that, for bounded solutions, $k(\d)$ is
always finite. A contradictions is obtained while assuming that
$t_j(\d) \to \bar t<1$ as $j \to \infty$ for arbitrarily small
$\d>0$ meaning a kind of ``complete blow-up" that is impossible
for conservation laws obeying the Maximum Principle.


This gives us a {\em global $\d$-deformation} in $\re \times
[0,1]$ of the solution $u(x,t)$, which is a discontinuous orbit
denoted by
 \beq
 \label{p3}
 \mbox{$
u^\d(x,t)= \{u_{j\d}(x,t) \,\,\, \mbox{for} \,\,\, t \in
[t_{j-1}(\d),t_j(\d)), \quad j=1,2,...,k(\d)\}.
 $}
 \eeq
Recall that the whole orbit (\ref{p3}) has been constructed by the
algebraic characteristic calculus using (\ref{p1}) only. Finally,
by an arbitrary {\em smooth $\d$-deformation}, we will  mean the
function (\ref{p3}) constructed by any sufficiently refined finite
partition $\{t_j(\d)\}$ of $[0,1]$,  without reaching a shock of
$S_-$-type at some or all intermediate points $t=t_{j}^-(\d)$.


We next say that, given a solution $u(x,t)$, it is {\em stable
relative smooth deformations}, or simply $\d$-stable ({\em
$\d$eformation-stable}), if for any $\e>0$, there exists
$\d=\d(\e)>0$ such that, for any finite $\d$-deformation of $u$
given by (\ref{p3}),
 \beq
 \label{p4}
 \mbox{$
\iint |u-u^\d| < \e.
 $}
 \eeq
 Then we have the following simple observation:

\begin{proposition}
 \label{Pr.D}
 Let under given hypothesis, a weak solution $u(x,t)$ of the Cauchy
 problem $(\ref{3})$, $(\ref{2})$ be $\d$-stable. Then it is entropy.
  \end{proposition}

Indeed, if $u(x,t)$ is not entropy, then there exists $t_* \in
(0,1]$ such that $u(x,t_*)$ does not satisfy (\ref{Ent2}), i.e.,
this profile has a finite non-entropy shock of the type $S_+$ at
some point $x_* \in \re$. Since those shocks cannot be reproduced
with arbitrary accuracy $\e$ in $L^1$ by the characteristic system
(\ref{p1}), any $\d$-deformation $u^\d$ at $t=t_*$ must stay
$\e_0>0$ away from $u(x,t_*)$ for arbitrarily small $\d>0$.

Of course, this construction does not play a role for conservation
laws with well-developed entropy theory, which establishes
existence of a semigroup of $L^1$-contractions of entropy
solutions. Obviously, this strong contractivity property
guarantees also uniqueness of $\d$-entropy solutions.
 The situation is different for the NDEs:


\section{On $\d$-entropy  solutions (a test) of the NDE}
\label{Sect4}

 Thus, we are going to
develop and discuss some aspects of entropy solutions for
(\ref{1}).
 {\em without} using the idea of vanishing, $\e \to 0$, viscosity as
in  \cite[\S~7]{GalEng}
 \beq
 \label{ee1}
  u_t=(uu_x)_{xx}- \e u_{xxxx} \quad \mbox{in}
 \quad \re \times \re_+ \quad (\e>0).
  \eeq
A direct verification that the $\e$-approximation (\ref{ee1})
yields as $\e \to 0$ the correct Kruzhkov's-type entropy solution
leads to difficult open problems.
 We begin with:


\subsection{(${\bf 2m}$+1)th-order NDEs for any $m \ge 1$ DO NOT generate a
semigroup of contractions in $L^1$}

 A first naive
approach would be to try to create a standard entropy condition
for the NDE of, say, the following form
(cf.
 (\ref{KrIn})):
\beq
 \label{KrInN}
 \mbox{$
 |u-k|_t - \frac 12 \, \big[{\rm sign}(u-k)  (u^2-k^2)\big]_{xxx} \le 0 \quad \mbox{in} \,\,\,
 {\mathcal D}'(Q_1) \quad \mbox{for any} \,\,\, k \in \re.
  $}
  \eeq
Then  Kruzhkov's-type computations with (\ref{1}) are supposed to
be performed by using his fundamental idea of doubling the space
dimension; see a clear presentation in \cite[p.~402]{Tay}, with
some obvious adaptations of test functions involved.

One should avoid doing this bearing in mind that this approach
must end up with the contractivity property (\ref{sem1}), which
{\em cannot be true for any PDE of order larger than two,} since
these are associated with  manipulations based on the  Maximum
Principle for first-order or, at most, for second-order parabolic
PDEs. This means that semigroups of contractions in $L^1$ are not
available for such NDEs (\ref{N5}) with any $m \ge 1$.

\subsection{On smooth solutions and odd-order operator theory}
\label{S6.3}

Thus, we return to the Cauchy problem for the NDE (\ref{1}).
As we have mentioned,  unlike the first-order case (\ref{cle1}),
applying the $\e$-approximation as in (\ref{ee1}) leads to a
number of principal difficult problems and, in the maximal
generality (excluding special cases), does not give neither
existence of a solution via the family $\{u_\e\}$ nor uniqueness
of an ``$\e$-entropy" solution, \cite{GalEng}.

 We will
develop other concepts of solutions by different types of
approximations, and then the concept of uniqueness will be
attached to the nature of existence results.

\ssk

\noi\underline{\em On local semigroup of smooth solutions}.
Beforehand, it is of importance that, as the similarity solutions
 in \cite[\S~3]{3NDEI}
showed, the NDE (\ref{1}) does not admit a global in time solution
for any bounded $L^1$ data. This is in striking difference with
the conservation laws (\ref{3}), where such existence is
guaranteed by the Maximum Principle. Therefore, we restrict our
attention to weak solutions $u(x,t)$ in $Q_1$, where
 \beq
 \label{u01}
 u_0(x) \in C_0^\infty(\re) \quad \mbox{is sufficiently small.}
  \eeq

Then,  as the first step of a similar construction,  we have to
check that
 for such smooth initial data $u_0$,
there exists a unique local classical $C^{3,1}_{x,t}$ solution
$u(x,t)$ of (\ref{1}).
 Recall that characteristic methods similar to that in (\ref{p1})
 are not available for higher-order PDEs.
 This just means that (\ref{1}) generates a
standard {\em local} semigroup in the class smooth functions.
 These results are known for non-degenerate NDEs such as
 (\ref{Lin12}), and moreover the solutions are $C^\infty$ locally in time,
  \cite{Cr90, Cr92, Lev01}.
   Actually, these smoothing results can
  be viewed in conjunction with
  classic methods of analytic semigroups in PDE theory;
 see
 \cite{PG} and references in a more recent paper \cite{Esch06};
 see below.

\ssk

\noi\underline{\em Uniqueness and continuous dependence: an
illustration}.  Actually, in our construction, we will need just a
local semigroup of smooth solutions that is continuous is
$L^1_{\rm loc}$. The fact that this is generated by third-order
(or other odd-order NDEs) is illustrated by the following easy
example. Consider, for definiteness, the NDE
 \beq
 \label{Gmm1GG}
 u_t={\bf A}(u) \equiv u u_{xxx}, \quad u(x,0) = u_0(x) \in H^7(\re),
  \eeq
  where, without loss of generality, we take into account
  the principal higher-order term only.  According to the above results, we assume
  that $u(x,t)$ satisfies
 \beq
 \label{mm2}
  \mbox{$
  \frac 1C \le u \le C \quad (C >1)
   $}
   \eeq
   and is
  sufficiently smooth,
   $u \in L^\infty([0,T], H^7(\re))$ and $u_t \in L^\infty([0,T], H^4(\re))$.
   See details on such uniqueness results  in
   \cite[\S~3]{Cr92}.

Thus, assuming that there exists the second smooth solution
$v(x,t)$, we subtract the equations and obtain for the difference
$w=u-v$ the following:
 \beq
 \label{uw11}
 w_t=u w_{xxx} + v_{xxx} w.
  \eeq
  We next divide (\ref{uw11}) by $u \ge \frac 1C>0$, multiply by $w$ in
  $L^2$, so, after integrating by parts, 
  \beq
  \label{mm3}
  \mbox{$
  \int \frac{w w_t}u \equiv \frac 12\, \frac{\mathrm d}{{\mathrm d}t}
   \int \frac {w^2}u+ \frac 12 \, \int \frac{u_t w^2}{u^2}
=
   \int \frac{v_{xxx}w^2}u.
   $}
   \eeq
Therefore, using  the assumed regularity yields
 \beq
 \label{mm4}
 \mbox{$
 \frac 12\, \frac{\mathrm d}{{\mathrm d}t}
   \int \frac {w^2}u= \int \big(\frac{v_{xxx}}u - \frac 12 \, \frac{u_t}{u^2} \big)  \, w^2
 \equiv \int \big(\frac{v_{xxx}}u - \frac 12 \, \frac{u_{xxx}}{u} \big)  \, w^2
 \le C_1 \int \frac {w^2}u,
  $}
  \eeq
  where we use the fact that $u_{xxx}(\cdot,t)$, $v_{xxx}(\cdot,t)
  \in
  L^\infty([0,T])$.
  By Gronwall's inequality, (\ref{mm4}) implies that $w(t) \equiv 0$.
As usual, this construction can be translated to the continuous
dependence result in $L^2$ and hence in $L^1_{\rm loc}$.

\ssk

\noi\underline{\em On degenerate NDEs}. For degenerate NDEs such
as (\ref{1}) and for solutions of changing sign,
 the unique local smooth solvability is a technical result, which we do not
completely concentrate upon, and present below some rather formal
comments justifying such a local continuation. One of the main
difficulties of this local analysis, is that (\ref{1}) admits
solutions with finite interfaces and free boundaries, which
represent ``weak shocks" with quite tricky (smooth enough but not
$C^3_x$) behaviour.

Thus, in addition, except the shock waves, which we are mostly
interested in, the NDE (\ref{1}) is degenerate at $\{u=0\}$, so
that the local existence of sufficiently smooth solution must
include the demand of ``transversality" of all the zeros (a finite
number) of initial data $u_0(x)$ (or $u(x,t_j(\d))$ later on).
Here the transversality of the zero at, say, $x=0$ has a standard
meaning:
 $$
 u_0'(0) \not = 0.
 $$
For instance, for key applications, we may assume that $u_0(x)$ is
anti-symmetric, so $u(-x,t) \equiv - u(x,t)$, and hence the only
transversal zero is fixed at the origin $x=0$ only, i.e.
 \beq
 \label{uu12}
 u(0,t) \equiv 0, \quad \mbox{and}
 \quad u(x,t)>0 \,\,\, \mbox{for} \,\,\, x<0.
  \eeq
 Then, according to regularity results for odd-order PDEs
 \cite{Cr90, Cr92, Hos99, Miz06} (cf. \cite{PG, Esch06, Lun}),
 the linearization about sufficiently smooth $u_0(x)$ yields that
the possibility of local smooth extension of solution is governed
by the good spectral properties of the third-order linear operator
with the principal part
 \beq
 \label{pr1}
  \mbox{$
 {\bf P}_3^1= x \frac{{\mathrm d}^3}{{\mathrm d}x^3} \quad \mbox{for} \quad x \approx
 0^+.
   $}
   \eeq
 This type of degeneracy is not sufficient to destroy good spectral
properties of ${\bf P}_3^1$ that still will admit a discrete
spectrum and a compact resolvent in the corresponding weighted
space $\sim L^2_{1/x}$ for $x>0$. Note that the singular point
$x=0$ starts to generate a continuous spectrum for the operator
 \beq
 \label{pr2}
  \mbox{$
{\bf P}_3^n= x^n \frac{{\mathrm d}^3}{{\mathrm d}x^3} \quad (x>0)
 $}
 \eeq
  in the parameter range $n \ge 3$
only, i.e., for much stronger degeneracy than in (\ref{pr1}).
 Indeed, then the change $z=x^\a$ with $\a= \frac{3-n}3>0$
 transforms (\ref{pr2}) into the regular operator with the constant
 principal part
 \beq
 \label{ff1}
 {\bf P}_3= D_z^3 \quad \mbox{for} \quad z \approx 0^+,
 \eeq
 for which all necessary spectral properties are obviously valid,
 \cite{NaiPartI}.
 The finite interface
behaviour will be shown to correspond to $n=2$, so it is still in
the good range. Our conclusions here are based on the well-known
fact that the linear PDE
 \beq
 \label{2.6}
 u_t=u_{xxx}
  \eeq
 generates a smooth
  (analytic in a properly weighted $L^2$-space)
semiflow given by
  \beq
  \label{bbb1}
 u(x,t)= b(x-\cdot,t)*u_0(\cdot),
  \eeq
 where $b(x,t)$ is the fundamental solution
 \beq
  \label{bb123}
  \mbox{$
b(x,t)= t^{-\frac 13} {F}\big( x/{t^{\frac 13}}\big),
\quad\mbox{where} \,\,\, F={\rm Ai}(z), \quad F'' + \frac 13 \,F
\,z=0, \quad \int F=1.
 $}
 \eeq

Thus, for the  degenerate NDE (\ref{1}), the notion of
``sufficiently smooth solutions" should also include the
assumption of transversality, i.e., of local behaviour near zeros.
Of course, this is not that essential hypothesis that has a local
character, and, for instance, completely disappears for the
related non-degenerate NDEs such as (\ref{non1}),
  which also admits shocks and needs proper entropy theory
  (to be treated also).

\ssk

\noi\underline{\em On odd-order ordinary differential operators}.
In the above analysis, we need a detailed spectral theory of
third-order (or more generally, odd-) operators such as
 \beq
 \label{bb1S}
 {\bf P}_3= a(z) D_z^3 + b(z) D_z^2 + c(z) D_z + d(z)I, \quad z \in
 (-L,L) \quad (a(z) \ge c>0),
  \eeq
  with bounded coefficients. This theory is available in Naimark's classic
   book \cite[Ch.~2]{NaiPartI}. It was shown that for {\em regular
  boundary conditions} (e.g., for periodic ones that are regular for any order and that suit us
  well), operators admit a
  discrete spectrum $\{\l_k\}$, where the eigenvalues $\l_k$ are
  all simple for $k \gg 1$, and a complete in $L^2$ subset of
  eigenfunctions $\{\psi_k\}$ that create a Riesz  basis\footnote{This is
  G.M.~Kessel'man's (1964) and V.P.~Mikhailov's
(1962)
  result.}.
  This makes it possible to use standard eigenfunction expansion
  techniques;
  see necessary  details and references at the end of Ch.~2 therein.

The eigenvalues of (\ref{bb1S})
 have the asymptotics
 \beq
 \label{bb2}
 \l_k \sim (\pm 2\pi k {\rm i})^3 \quad \mbox{for all \,\, $k \gg 1$}.
 \eeq
In particular, this means that ${\bf P}_3 - aI$ for any $a \gg 1$
is not a sectorial operator that makes suspicious referring to the
analogies with analytic theory \cite{PG, Esch06, Lun} that is
natural for even-order parabolic flows.

Nevertheless, recall that (\ref{bbb1}) guarantees analyticity of
solutions that is now associated with the Airy-type operator
 \beq
 \label{BB1}
  \mbox{$
 {\bf B}_3= D_z^3 + \frac 13 \, z D_z + \frac 13 \, I
 \quad \mbox{in}
 \quad L^2_\rho(\re), \quad \rho(z)={\mathrm e}^{a|z|^{3/2}},
  $}
 \eeq
 where $a>0$ is sufficiently small; cf. a ``parabolic" version of
 such a spectral theory in \cite{Eg4}. It turns out that
 (\ref{BB1}) has the real spectrum (see \cite[\S~9]{2mSturm})
  $$
   \mbox{$
   \s({\bf B})=\big\{- \frac
 l3, \, l=0,1,2,...\big\},
  $}
  $$
   so that ${\bf B}-a I$ is sectorial for
 $a \ge 0$ ($\l_0=0$ is simple), and this justifies the fact that
 (\ref{bbb1}) is an analytic flow.

Note also that analytic smoothing effects are known for
higher-order dispersive equations with operators of principal
type, \cite{Tak06}.
  This suggests to treat
 (\ref{Gmm1GG}) by classic approach as in Da
 Prato--Grisvard \cite{PG} by linearizing about a sufficiently
 smooth $u_0=u(t_0)$, $t_0 \ge 0$, by setting $u(t)=u_0+v(t)$ giving the
 linearized equation
  \beq
  \label{leq1}
  v_t= {\bf A}'(u_0)v + {\bf A}(u_0)+ g(v), \quad t > t_0; \quad
  v(t_0)=0,
   \eeq
   where $g(v)$ is a quadratic perturbation. Using good semigroup properties of
    ${\mathrm e}^{{\bf A}'(u_0)t}$, this makes it possible to
    study local regularity properties of the
    integral equation
     \beq
     \label{leq2}
     \mbox{$
     v(t) =\int\limits_{t_0}^t{\mathrm e}^{{\bf A}'(u_0)(t-s)}({\bf A}(u_0)+
     g(v(s)))\, {\mathrm d}s.
      $}
       \eeq
It is key that the necessary smoothness of solutions  demands the
fast exponential decay of solutions $v(x,t)$ as $x \to \infty$,
since one needs that $v(\cdot,t) \in L^2_\rho$; cf. \cite{Lev01},
where $C^\infty$-smoothing also needs an exponential-like decay.
Equations such as
(\ref{leq2}) can be used to guarantee local existence of smooth
solutions of a wide class of odd-order NDEs.


 Thus, we state the following conclusion to be used later on:
 \beq
 \label{cc11}
  \begin{matrix}
 \mbox{any sufficiently smooth solution $u(x,t)$ of (\ref{Gmm1GG}),
 (\ref{mm2}) at $t=t_0$,}\ssk\\
 \mbox{can be uniquely extended to some interval
 $t \in(t_0,t_0+\d)$,  $\d>0$.}
  \end{matrix}
  \eeq




\subsection{Global solutions by Galerkin method}

Here we demonstrate the application of another classic approach to
nonlinear problems that, suddenly, in the present case of unclear
entropy nature of solutions of NDEs and the open uniqueness
problem, gives a partial answer to both. We mean the Galerkin
method that was the most widely used approach for constructing
weak solutions via finite-dimensional approximations; see Lions
\cite{LIO} with many applications therein.

 Thus, by this classic
theory of nonlinear problems, under the assumption (\ref{u01}) and
others, if necessary, let us perform a standard construction of a
compactly supported (for simplicity)
    solution by
Galerkin method using the basis $\{\psi_k\}$ of eigenfunctions of
the regular linear operator ${\bf P}_2=D_x^2<0$ with the Dirichlet
boundary conditions,
 \beq
 \label{m10}
 \psi''= \l_k \psi, \quad
 \psi=0 \,\,\, \mbox{at} \,\,\, x=\pm L \quad \Longrightarrow \quad
 \l_k \sim -k^2.
  \quad \mbox{and}
  \quad  u=0 \,\,\, \mbox{at} \,\,\, x=L.
   \eeq
As an alternative, it is
 curious that, for our purposes,
possible (and more convenient for some reasons) to use the
eigenfunction set
 of the operator ${\bf P}_4=-D_x^4<0$ again  with the
  Dirichlet conditions
   $$
   \mbox{$\psi=\psi_x=0$ \,\,at\,\, $x=\pm L$}.
   $$
 Special Galerkin bases associated with higher-order operators
 ${\bf P}_6=D_x^6<0$ are also may be convenient; see applications
 to third-order linear dispersion equations in \cite{Lar06}.

 In all these self-adjoint  cases, the
eigenfunctions form a complete and closed set in $L^2$;
 see classic theory of ordinary differential
operators in Naimark \cite[p.~89]{NaiPartI}.



On the other hand,  looking more natural  choice of the
third-order operator ${\bf P}_3=D_x^3$ for Galerkin approximation
of (\ref{1}) will cause a difficult problem, since for the
third-order PDE with the principal operator as in (\ref{1}),
 \beq
 \label{ss1N}
 u_t = a(x,t) u_{xxx}+... \quad (a \mapsto u)
 \eeq
 with $a>0$,
   proper setting for the
 IBV problem includes the Dirichlet
  conditions (see Faminskii \cite{Fam02} for details and a survey)
  \beq
  \label{m2}
  u=u_x=0 \,\,\, \mbox{at} \,\,\, x=-L \quad \mbox{and}
  \quad  u=0 \,\,\, \mbox{at} \,\,\, x=L.
   \eeq
For $a<0$, the boundary conditions must be swapped, so that the
proper setting of the problem depends on the unknown sign of
solutions. Here, the fact that ${\bf P}_3=D_x^3$ is not
self-adjoint is not essential since, relative to adjoint basis
$\{\psi_k^*\}$,   the closure and completeness of the
bi-orthonormal generalized eigenfunction sets remain valid.

Actually, the choice of  linear operators ${\bf P}_2=D_x^2$, ${\bf
P}_4=-D_x^4$, or others, is not of principal importance if we are
looking for compactly supported solutions
  \beq
 \label{m5}
 u \in C_0^\infty((-L,L) \times [0,1]).
  \eeq
It should be noted that the control of finite propagation property
in (\ref{1}) is difficult and is an essential part of our further
analysis.
 For instance, we also can fix periodic boundary conditions that
 are always  {regular}, \cite[Ch.~2]{NaiPartI} (it is curious that (\ref{m2}) are not).

Thus, we construct a sequence $\{u_m\}$ of approximating Galerkin
solutions of (\ref{1}), (\ref{2}) in the form of finite sums
 \beq
 \label{f1}
  \mbox{$
 u_m(x,t)= \sum\limits_{k=1}^m C_k(t) \psi_k(x),
  $}
   \eeq
 where $\{C_j\}$ solve the quadratic dynamical systems
 \beq
 \label{f2}
  \mbox{$
C_j'= \sum\limits_{(k,l)} C_k C_l J_{klj}, \quad \mbox{where}
\quad J_{klj}= \langle \psi_k \psi_l', \psi_j'' \rangle = \l_j
 \langle \psi_k \psi_l', \psi_j \rangle.
 $}
 \eeq
  For
 the conservation law (\ref{3}), the DS takes the same form
 as in (\ref{f2}), with the only difference that
  \beq
  \label{bb1}
J_{klj}= -\langle \psi_k \psi_l', \psi_j \rangle.
 \eeq

The identity  (\ref{111}) for $u_m$ takes the form
 \beq
  \label{ss1NNN}
  \mbox{$
  \sum\limits_{(k)} \frac 1{|\l_k|} \, C_k^2(t)= c_{0m}=
  \sum\limits_{(k)} \frac 1{|\l_k|} \, C_k^2(0), \quad t>0.
 $}
 \eeq
This guarantees global existence of the solutions $u_m(x,t)$
showing that
 \beq
 \label{ss2N}
 C_k(t) \quad \mbox{do not blow-up and exist for all $t>0$.}
  \eeq

Since $\psi_k$ are given by $\sin(\l_k x)$ or $\cos(\l_k x)$, a
lot of coefficients $J_{klj}$ vanish. For instance, if $u_0(x)$ is
odd,  we take all the sin-functions,
 $$
  \mbox{$
 \psi_k(x)= \frac 1{\sqrt L} \, \sin\big(\frac{k \pi x}L\big), \quad
 \mbox{with} \quad \l_k= - \frac{k^2 \pi^2}{L^2}, \quad k=1,2,...
 \, .
  $}
  $$
The non-zero coefficients $J_{klj}$ occur iff
 $
 k=j$,
 $ l=2j,
  $
  where (\ref{f2}) becomes  simpler,
   \beq
   \label{lp1}
    \mbox{$
   C_j'= \frac{2\pi^3 j^3}{L^{7/2}} \, C_j C_{2j}, \quad j \ge 1.
     $}
     \eeq
It is curious that (\ref{lp1}) yields the following feature of a
``maximum principle":
 \beq
 \label{lp2}
 \mbox{$
{\rm sign} \, C_j(t) = {\rm sign} \, C_j(0), \quad j \ge 1.
 $}
 \eeq


Other {\em a priori} estimates
 are obtained  by multiplying
(\ref{1}) in $L^2$ by $u$ and $u_{xx}$ yielding the identities
 \beq
 \label{m3}
  \mbox{$
\frac 12 \, \frac{\mathrm d}{{\mathrm d}t}\,  \int u^2= -\frac 12
\, \int (u_x)^3, \quad
 \frac 12 \, \frac{\mathrm d}{{\mathrm d}t}\,  \int (u_x)^2=- \frac 52
\, \int u_x(u_{xx})^2. $}
 \eeq
Then some  interpolations of various terms in the identities
(\ref{m3}) are necessary.

Thus, the sequence of ``regularized" solutions (Galerkin
approximations) $\{u_m(x,t)\}$ is globally defined, and
 \beq
 \label{g1}
 \{u_m\} \quad \mbox{is uniformly bounded in
 $L^\infty([0,1];H^{-1})$}.
  \eeq
Therefore, along a subsequence, $\{u_m\}$ converges to $u$
weakly-* in  $L^\infty([0,1];H^{-1})$, and, in addition, strongly
in
 $H^{-1}([0,1];H^{-2})$, in view of compact embedding.
This gives a weak solution. As usual, the better regularity comes
from the special choice of Galerkin's basis employed. We do not
stress attention to this (bearing in mind local
$C^\infty$-smoothing for non-degenerate NDEs). See \cite{Lar06}
for rather exotic Galerkin bases  applied
 to KdV type equations. Recall
that, globally, smoothing is not available, since this
construction is specially oriented to include shocks of
$S_-$-type.

\ssk

 \noi{\bf Remark 1.}
 Obviously, the estimate (\ref{ss1NNN}) does not and cannot
prevent gradient catastrophe, which means that
 \beq
 \label{kk1}
 \mbox{$
\|u_x(t)\|_2^2 = \sum |\l_k| C_k^2(t) \to + \infty
 \quad \mbox{as} \quad t \to T^- \le 1.
 $}
 \eeq
Notice that for (\ref{1}) there is an opportunity to create
blow-up of the solutions $u(\cdot,t)$ itself (possibly together
with (\ref{kk1})), where
 \beq
 \label{kk1NN}
 \mbox{$
\|u(t)\|_2^2 = \sum  C_k^2(t) \to + \infty
 \quad \mbox{as} \quad t \to T^- \le 1.
 $}
 \eeq
This does not happen if a finite shock appears via the
self-similar patterns such as \cite{3NDEI}
 \beq
 \label{2.1}
 u_-(x,t)=g(z), \,\, z= x/(-t)^{\frac 13},
  \mbox{$
 \quad (g g')''= \frac 13 \, g'z,
  \,\,\,
  f(\mp
  \infty)=\pm 1.
   $}
   \eeq
 Indeed,  by the first
identity in (\ref{m3}), there appears an integrable singularity,
 \beq
 \label{m3ss}
\mbox{$
 \frac 12 \, \frac{\mathrm d}{{\mathrm d}t}\,  \int\limits_{-1}^1 u^2
 \sim
-\frac 12 \, \int\limits_{-1}^1 (u_x)^3 \sim (-t)^{- \frac 23}
\int\limits_{-(-t)^{-1/3}}^{(-t)^{-1/3}} (g')^3\, {\mathrm d}z \in
L^1((-1,0)),
 $}
 \eeq
so that $\|u(0^-)\|_2^2$ remains finite. Here in (\ref{m3ss}) one
needs to use the asymptotics of the Airy function
\cite[\S~3]{3NDEI}, so that the integral therein diverges but  its
rate,
 $$
 \mbox{$
\Big| \int\limits_{-(-t)^{-1/3}}^{(-t)^{-1/3}} (g')^3 \, {\mathrm
d}z \Big| \sim O((-t)^{-\frac 1{10}}),
 $}
 $$
 is sufficient for the integrability.

\ssk

\noi{\bf Remark 2.} Using the dynamical system (\ref{f2}) instead
of the NDE (\ref{1}) suggests to develop a formal calculus of the
corresponding sequences, where, on identification,
 \beq
 \label{ss1P}
 \mbox{$
u= \sum C_k \psi_k  \quad \Longrightarrow \quad u=\{C_k\}
 $}
 \eeq
belongs to the little Hilbert space $h^{-1}_{\bf P}$ with the
metric
 \beq
 \label{ss2P}
 \mbox{$
\|u\|_{\bf P}^2= \sum \frac 1{|\l_k|}C_k^2.
 $}
 \eeq
Then (\ref{ss1NNN}) guarantees that
  \beq
 \label{ss3P}
 \mbox{$
u(t) \in h^{-1}_{\bf P} \quad \mbox{for all $t \ge 0$},
 $}
 \eeq
meaning global solvability. Moreover,  the embedding $h^{-2}_{\bf
P} \subset h^{-1}_{\bf P}$ is compact since $\frac 1{|\l_k|} \sim
\frac 1{k^2}$ \cite{LustSob} (for $h^{-2}_{\bf P}$, the metric
contains $ \frac 1{|\l_k|^2}$ in (\ref{ss2P})), so that  we can
use the same Galerkin approximation method to construct suitable
solutions. In this space, the blow-up formation of shocks means
(\ref{kk1}).

\ssk

\noi{\bf Remark 3.} Writing the $N$-dimensional NDE (\ref{ndeN})
for compactly supported  $u_0$ as
 \beq
 \label{ndeN1}
  \mbox{$
 (-\D)^{-1} u_t = - \frac 12\,
 \frac{\partial}{\partial x_1}\, u^2
 $}
 \eeq
 with the standard definition of the linear operator $(-\D)^{-1}$
 in $L^2(\O)$, $\O$ is sufficiently large, and multiplying
 (\ref{ndeN1}) by $u$ yields
 the same conservation identity (\ref{111}).
 Some  concepts
developed above  can be also adapted to the equations in $\ren$,
though shock wave formation phenomena become more involved and are
in general unknown.

\subsection{$\d$-entropy solutions (a test) for the NDE}


 Assuming that the local smooth
solvability problem above is well-posed, we now present the
corresponding definition that will be
 applied to  particular weak solutions. Recall that the topology
 of convergence, $L^1_{\rm loc}$ at present, for (\ref{1}) was justified
by a similarity analysis presented in \cite[Prop.~3.2]{3NDEI}.
 For
other NDEs, the topology may be different that can be
 a difficult
problem.

  \ssk

 \noi{\bf Definition \ref{Sect4}.1.} A weak solution $u(x,t)$ of
 the Cauchy problem $(\ref{1})$, $(\ref{2})$ is
 called $\d$-entropy if there exists a sequence of its smooth $\d$-deformations
 $\{u^{\d_k}, \, k=1,2,...\}$, where $\d_k \to 0$, which converges in $L^1_{\rm loc}$ to
 $u$ as $k \to \infty$.

 \ssk

Note that this is slightly weaker (but equivalent) to the
condition of $\d$-stability. The construction of global
$\d$-deformation of $u$ is performed along the lines of (i)--(k)
in Section \ref{S6.2}. The only difference is that local
$\d$-deformations can lead to complete blow-up for the NDE
(\ref{1}), as explained in
 \cite[\S~4.2]{3NDEI}.
  To avoid this, one needs either to impose the
condition (\ref{u01}) or specially assume that complete blow-up
cannot occur under slight deformation of the data, or while
performing its $\d$-deformation  with any sufficiently small
$\d>0$. We call such solutions {\em $\d$-extensible} (the
definition assumes that $u$ is $\d$-extensible).

\ssk

\noi\underline{\em On $\d$-entropy test and uniqueness}.
 First of all, we again note that any uniqueness (and entropy) results for such
 NDEs are not acievable in principle, \cite{Gal3NDENew}.
 Therefore, we use the above results as a basis of the so-called
 "$\d$-Entropy Test" for testing shock an rarefaction waves; see
 first applications below.

\ssk

\noi \underline{\em $\d$-entropy solutions: motivation of the
term}. Let us explain why solutions are called $\d$-{\em entropy},
while
 we do not  use any evolution integro-differential inequality
 such as (\ref{KrIn}).
 It turns out that the NDE (\ref{1}) itself  contains the
 right evolution choice  of the admitted type shocks in the
 class of smooth solutions (precisely this makes sense of Definition
 \ref{Sect4}.1).

For instance, as a rough explanation, assume that at $x=0$ the
shock $S_+$ is going to appear at $t=1^-$ from a smooth solution
$u(x,t)$ such that $u(x,1^-)$ remains smooth everywhere except
$x=0$; e.g., for simplicity, we assume that
 \beq
 \label{33.1}
 u(x,1) \approx S_+(x) \quad \mbox{in a neighbourhood $x \in (-\d,\d)$},
  \eeq
  together with necessary derivatives $u_x$ and $u_{xx}$ that are
  assumed to be small at $x = \pm \d$. Here  $\d>0$ is also a small
  constant, so our illustration is of local nature. Multiplying
  (\ref{1}) by $u$ and integrating over $(-\d,\d)$ for $t \approx
  1^-$ yields the following main terms:
   \beq
   \label{h1}
   \mbox{$
   \frac 12 \,  \frac{\mathrm d}{{\mathrm d}t}\,
   \int (u^2-1)=  \frac 12 \int (u_x)^3 +... > 0
   \quad\big(\mbox{or $<0$ for $S_-(x)$ at $t=1$}\big),
   $}
    \eeq
    since $u_x$ must be essentially positive on profiles $u(x,t)$ that
     smoothly approximate
    $S_+(x)$. One can see that (\ref{h1}) evolutionary prohibits
    stabilization to $S_+(x)$ as $t \to 1^-$, when $u^2 \to 1$ in $L^1_{\rm loc}$.
     More rigorously
    \cite[\S~7.2]{GalEng}, the same
  negative result is established using the weaker topology of
  $H^{-1}$, where multiplication applies to the non-local equation
  (\ref{10}).

Similarly, we arrive at no contradiction while using (\ref{h1}) to
describe stabilization to $S_-(x)$, since then $u_x$ is
essentially negative. In fact,  (\ref{h1}) reflects a finite-time
formation of the singular shock $S_-$ (the gradient catastrophe)
for the NDE (\ref{1}) that was described  in \cite[\S~3]{3NDEI}
in greater detail.

Thus, using smooth deformations guarantees (via smoothness, that
is important) the preservation of the natural {\em local
entropies} such as inequalities like (\ref{h1}) and the opposite
one for $S_-$, so we call the constructed solutions $\d$-entropy.


\ssk

\noi\underline{\em First easy application  of $\d$-entropy test}.
As a first application, we have:

\begin{proposition}
\label{Pr.E} Shocks $S_-(x)$ and $H(-x)$ are $\d$-entropy.
 \end{proposition}

 The result follows from the properties of similarity solutions
 (\ref{2.1}), which, by  shifting the blow-up time $T \mapsto T+\d$,
  can be used  as their local smooth
 $\d$-deformations
 at any point $t \in [0,1)$. For $H(-x)$,
   we will need an extra
 approximation of similarity profile $g(z)$ with finite interface
 at some $z=z_0$, at which it is not $C^3$,  by sufficiently smooth profiles.


Let us use the negation in the following form:

 \ssk

 \noi{\bf Definition  \ref{Sect4}.2.} A weak solution $u(x,t)$ of
 the Cauchy problem $(\ref{1})$, $(\ref{2})$ is not $\d$-entropy
  if it is not $\d$-stable.

 \ssk

\begin{proposition}
\label{Pr.NE} Shocks $S_+(x)$ and $H(x)$ are not  $\d$-entropy.
 \end{proposition}

Indeed, taking initial data $S_+(x)$ and constructing its smooth
$\d$-deformation via the self-similar solution
\cite[\S~3.4]{3NDEI}
  \beq
  \label{2.14}
   \mbox{$
  u_+(x,t) = g(z), \,\,\,
   z=x/t^{\frac 13},
 \quad  (g g')''=- \frac 13 \, g'z,\,\,\,
 f(\mp
  \infty)=\mp 1.
   $}
   \eeq
Performing time-shifting $t \mapsto t+\d$, we obtain the global
$\d$-deformation $\{u^\d=u_+(x,t+\d)\}$ which goes away from
$S_+$.


Thus, we have shown that, at least, the idea of $\d$-deformations
allows  us to distinguish basic $\d$-entropy and non-entropy
shocks {\em without} any use of mathematical manipulations
associated with standard entropy inequalities, which are illusive
for higher-order NDEs (and nonexistent  in principle
\cite{Gal3NDENew}).
 \section{Compactons are $\d$-entropy solutions}
 \label{Sect6}

Without loss of generality, we treat this question for a
particular NDE. Namely,
 consider the following $K(2,2)$ equation:
 \beq
 \label{61}
 u_t = (uu_x)_{xx} + 4u u_x \quad \mbox{in}
 \quad \re \times \re_+.
  \eeq
Its compacton solution has the explicit form \cite{RosH93},
 \beq
 \label{62}
 u_{\rm c}(x,t)= f_{\rm c}(x+ 3 t), \quad \mbox{where}\quad
f_{\rm c}(y)= \left\{
 \begin{matrix} 2 \cos^2\big(\frac y 2\big) \,\,\,\, \mbox{for}
 \,\,\, |y| \le \pi, \ssk\\
 \quad 0 \quad \qquad \mbox{otherwise}.
  \end{matrix}
  \right.
  \eeq
  This is an example of a compactly supported weak solution of
 equation (\ref{61}).  One can see that
at the interface points $y= \pm \pi$, the profile $f_{\rm c}(y)$
is just $C^{1,1}_y$, i.e., the first derivative $f_{\rm c}'(y)$ is
Lipschitz. Therefore, it is not a classical $C^{3,1}_{x,t}$
solution of the PDE and has weak singularities at $y= \pm \pi$, so
one needs to check whether it is an entropy solution.
 In
addition, the ``flux" $(f f')'$ is continuous at those points,
though this does prove nothing.

We now use the concept of $\d$-entropy solutions from Sections
\ref{Sect40} and \ref{Sect4}.

\begin{proposition}
 \label{Pr.C}
The compacton $(\ref{62})$ is a $\d$-entropy solution of the NDE
$(\ref{61})$.
 \end{proposition}

 \noi{\em Proof.} We are going to show that there exists smooth
 $\d$-deformations of $u_{\rm c}$ for arbitrarily small $\d>0$. The general TW
solutions as in (\ref{62}) with $\l=-3$ yields the ODE
 \beq
 \label{63}
3 f'= (ff')''+2 (f^2)' \quad \Longrightarrow \quad 3f=(ff')'+ 2f^2
+ C_\d,
 \eeq
 where we chose the constant of integration to be
 \beq
 \label{64}
   C_\d= 3\d- 2\d^2>0.
    \eeq
 One can see on the phase plane in the variables $\{f^2,(f^2)'\}$
 that the ODE (\ref{63}), (\ref{64}) has a strictly positive and hence analytic
solution $f_\d$ satisfying
 \beq
 \label{65}
 f_\d(y) \to \d^+ \quad \mbox{as}
 \quad y \to \pm \infty \,\,\, \mbox{exponentially fast}, \quad \mbox{and}
 f_\d \to f_{\rm c} \quad \mbox{as \, $\d \to 0^+$}
 \eeq
  uniformly in
 $\re$.
  According to Definition \ref{Sect4}.1,  (\ref{65}) implies that
  $u_{\rm c}$ is an entropy
   solution of (\ref{61}), as well as $f_{\rm c}$ is  G-admissible for the third-order ODE
 in (\ref{63}). $\qed$



  It is crucial that Proposition \ref{Pr.C} justifies that the
  $K(2,2)$ equation (\ref{61}), the NDE (\ref{1}), and many others
  with similar principle degenerate third-order operators
  possess
  {\em finite propagation of interfaces} for entropy solutions.

 It is worth recalling again that, regardless the existence of such
nice smooth compactons (\ref{62}), the generic behaviour for the
RH equation (\ref{61}), for other data, includes formation of
shocks  in finite time, with the local similarity
mechanism as in \cite[\S~3.1]{3NDEI}. 

\section{On extensions to other related  NDEs}
 \label{Sect9}




\subsection{Shocks for the non-degenerate NDE}

\begin{figure}
\centering
\subfigure[shock profiles]{
\includegraphics[scale=0.52]{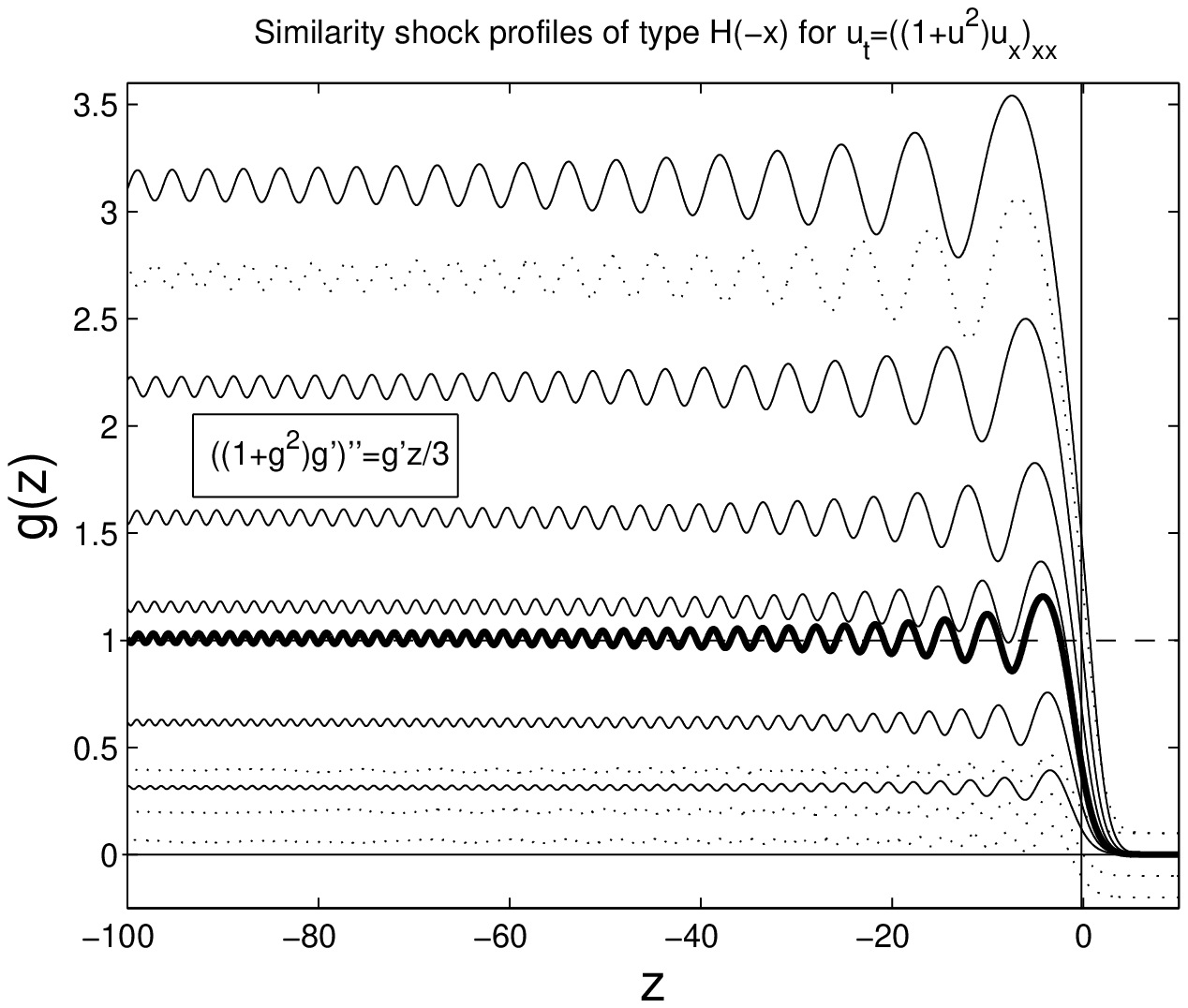} 
}
\subfigure[soliton]{
\includegraphics[scale=0.52]{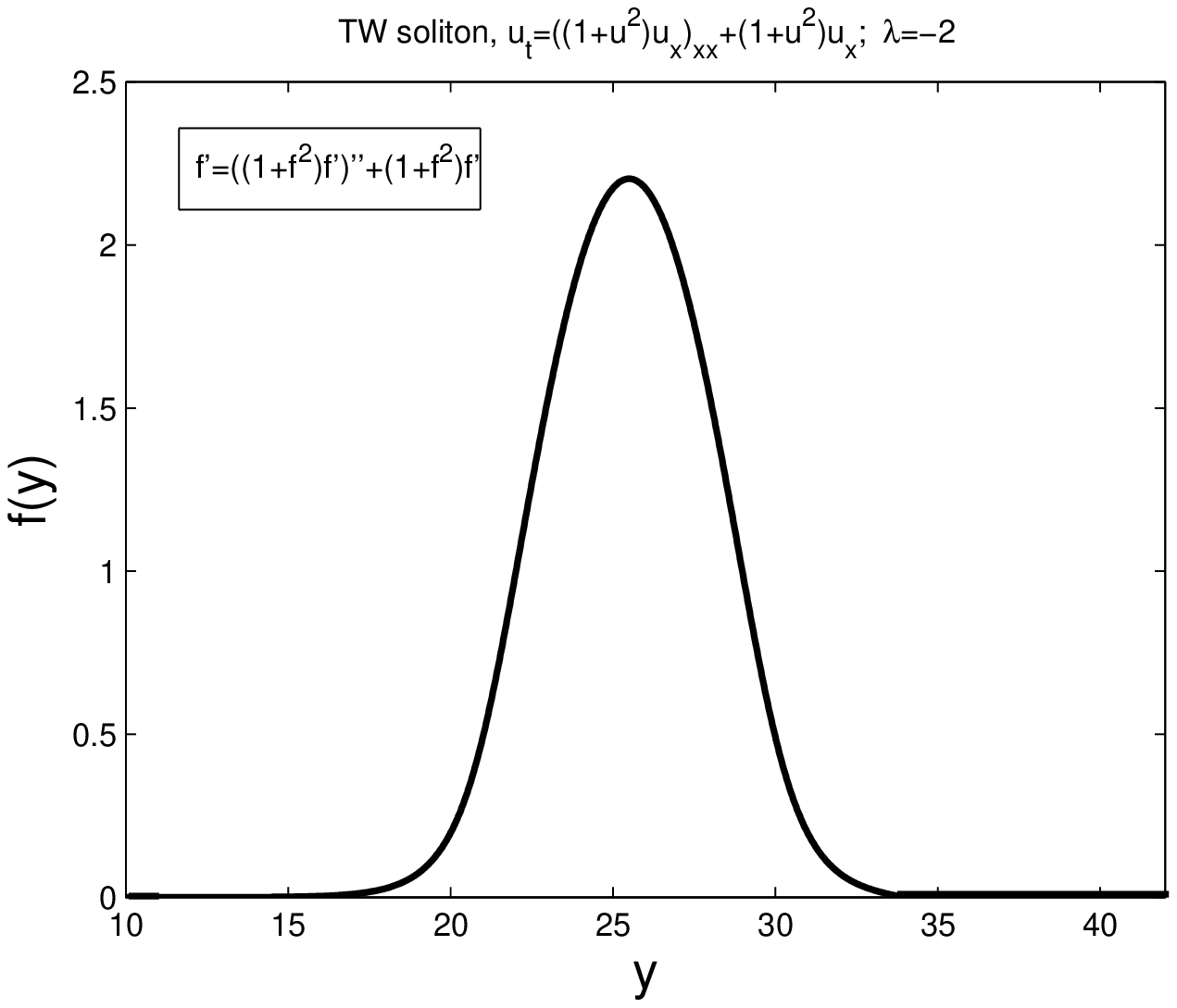} 
}
 \vskip -.4cm
\caption{\rm\small Similarity profiles for the NDE (\ref{non1}):
 shock  profiles satisfying the ODE in (\ref{gg1}) (a), and TW soliton
satisfying (\ref{non21})
 with $\l=-1$.}
 \vskip -.3cm
 \label{FQQ1}
\end{figure}

We begin with the simpler model (\ref{non1}) that appeared in
Section \ref{S6.3}  while we discussed the possibility of
extensions of sufficiently smooth solutions for defining
$\d$-deformations. Indeed, for (\ref{non1}), this is much easier.
On the other hand, obviously, as an NDE, this admits shocks via
standard similarity solutions
 \beq
 \label{gg1}
  \mbox{$
 u_-(x,t)=g(z), \quad z=x/(-t)^{\frac 13}
  \quad \Longrightarrow \quad
   ((1+g^2)g')''= \frac 13 \, g'z.
   $}
   \eeq
This ODE is studied as usual. Figure \ref{FQQ1}(a) shows a few
similarity profiles
 satisfying
  \beq
  \label{kkkk1}
  \mbox{$
 g(z) \sim z^{-\frac 14} {\mathrm e}^{-a_0 z^{3/2}}\to 0 \quad \mbox{as}
 \quad z \to + \infty, \quad \mbox{where} \,\,\, a_0= \frac 2{3 \sqrt 3},
 $}
 \eeq
that create as $t \to 0^-$ the shocks $\sim H(-x)$. By dotted
lines, we indicate there other profiles $g(z)$, for which
 $g(+\infty) \not = 0$.
 For the sake of comparison with compactons, in Figure
 \ref{FQQ1}(b),  we present the {\em soliton} of the related NDE
 \beq
 \label{non22}
 u_t=((1+u^2)u_x)_{xx} + (1+u^2)u_x, \quad \mbox{where}
  \eeq
   \beq
   \label{non21}
   u_{\rm s}(x,t)=f_{\rm s}(y), \quad y= x- \l t
   \quad \Longrightarrow \quad -\l f'=((1+f^2)f')''+ (1+f^2)f'.
    \eeq
    The soliton profiles
have now  exponential decay for $\l <-1$,
 $$
 f_{\rm s}(y) \sim {\mathrm e}^{-a_0|y|} \to 0 \quad \mbox{as}
 \quad |y| \to +\infty, \quad a_0=\sqrt{|1+\l|}.
  $$



 \subsection{
 $\d$-entropy approach to the NDE--(2,1)}

 For the non-fully divergent PDE (\ref{3,1})
   we  also apply  the    $\d$-entropy
   to prove existence and
uniqueness via suitable approximations.

\ssk

\noi\underline{\em On Galerkin method}. Constructing Galerkin
approximations, we face a new technical difficulty in passing to
the limit since a uniform estimate such as (\ref{ss1NNN}) is not
available for solutions \cite[\S~4]{3NDEI}
 \beq
 \label{3.1}
  \mbox{$
 u_\a(x,t)=(-t)^\a g(z), \quad z= x/(-t)^{\b}, \quad \b=
 \frac{1+\a}3 \quad(\a \in \re).
  $}
  \eeq
 Nevertheless, we can establish some extra estimates by using the
corresponding DS (\ref{f2}), where
 $
 J_{klj}=- \l_l \langle \psi_k \psi_l, \psi_j' \rangle.
  $
E.g., for odd data, the simpler system similar to (\ref{lp1}),
 $$
C_j'= \g_0 j^3 C_j C_{2j}, \quad 1 \le j \le m,
 $$
 implies that, for $m$ even,
  \beq
  \label{CC1}
   \mbox{$
  C_{\frac m 2}'= \g_0 \big(\frac m 2\big)^3 C_{\frac m 2} C_{m} \,\, \mbox{and} \,\, C_j'=0, \,\,\,
  j > \frac m 2 \,\,\Longrightarrow \,\,
  C_{\frac m 2}(t)= C_{\frac m 2}(0){\mathrm e}^{\g_0 (\frac m2)^3 C_{m}(0)\, t}.
   $}
  \eeq
Therefore, assuming that
 \beq
 \label{uu1}
  \mbox{$
 u_0 \in C_0^3 \quad \Longrightarrow \quad |C_m(0)| \le \frac
 {c_*}{m^3}, \,\,\, m \ge 1 \quad (c_* >0),
 $}
  \eeq
 we obtain from (\ref{CC1}) a uniform bound on the Galerkin
 coefficients $\{C_j\}$, and hence a local weak solution.

  \ssk



\begin{figure}
\centering
\subfigure[for $S_-(x)$]{
\includegraphics[scale=0.52]{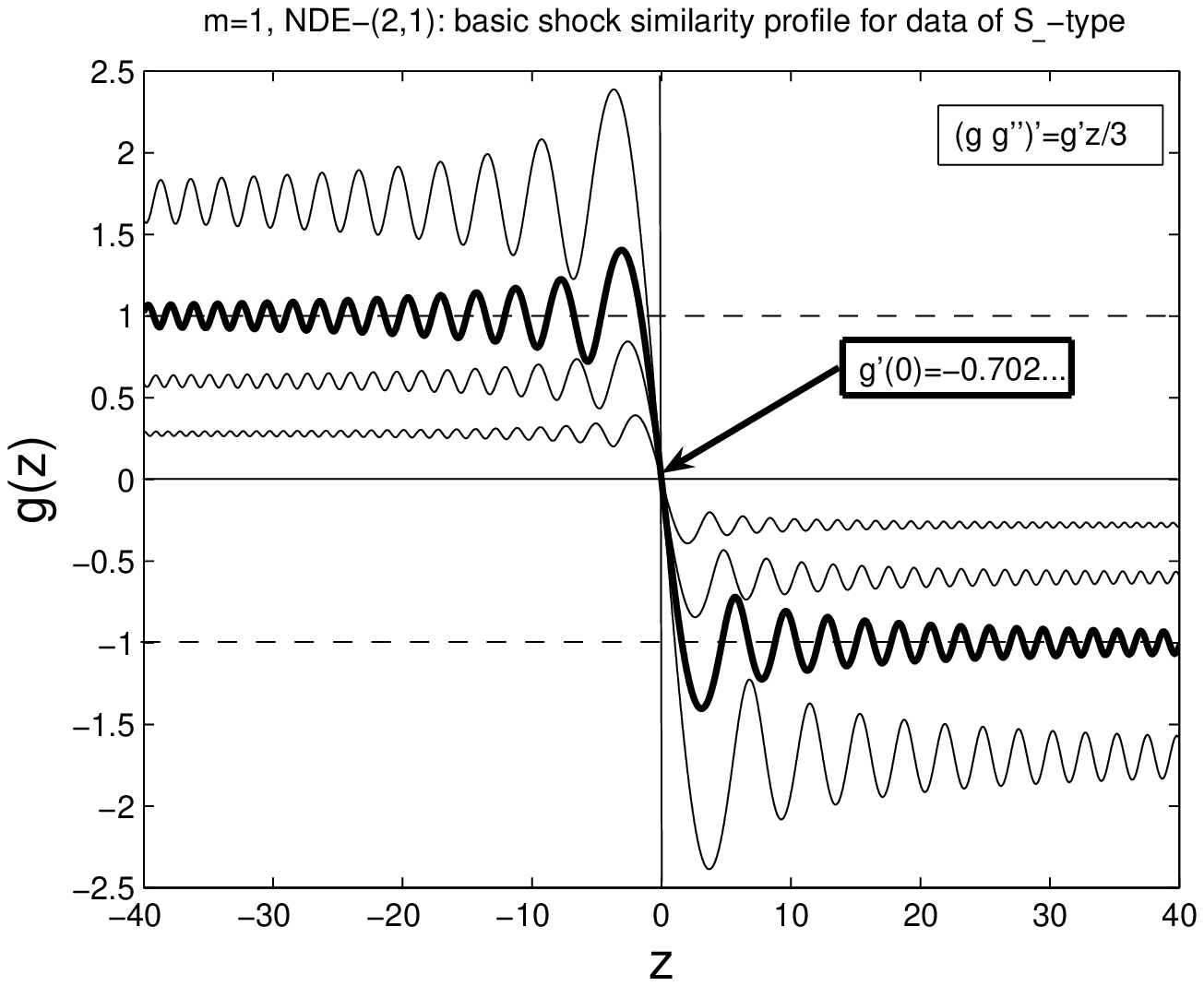} 
}
\subfigure[for $H(-x)$]{
\includegraphics[scale=0.52]{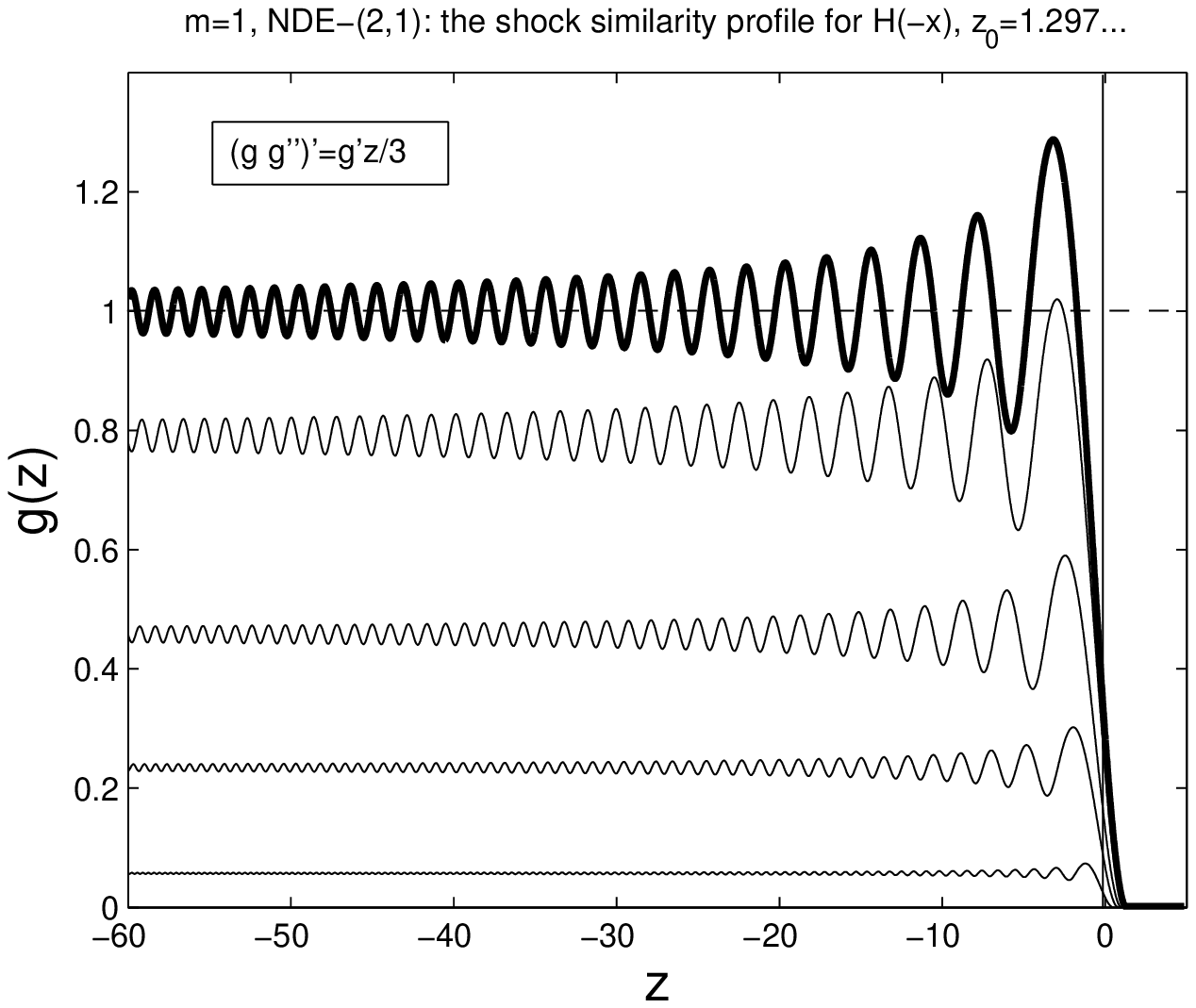} 
}
 \vskip -.4cm
\caption{\rm\small The ODE (\ref{661}):  the shock similarity
profiles including the unique solution (boldface line) for data
$S_-(x)$ (a), and  for $H(-x)$ with finite right-hand interface
(b).}
 \vskip -.3cm
 \label{F1N}
\end{figure}

\begin{figure}
\centering
\includegraphics[scale=0.70]{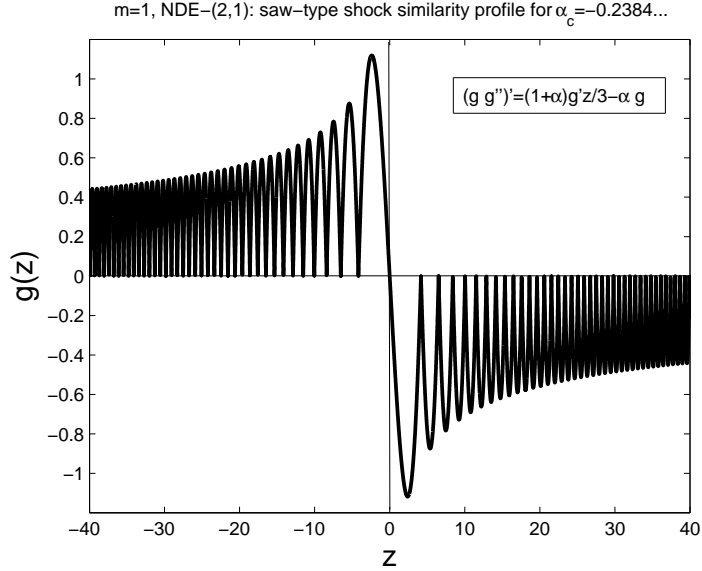}
 \vskip -.3cm
\caption{\small The saw-type similarity solution of the ODE in
(\ref{663}) for $\a_c=-0.2384...\,$.}
   \vskip -.3cm
\label{F2N}
\end{figure}

\noi\underline{\em Shocks and compactons exist}. On the other
hand, regardless its non-full divergence and nonexistence of any
obvious conservation laws, the NDE
(\ref{3,1}) allows a similar treatment of shocks and rarefaction
wave as for (\ref{1}).
For instance, formation of finite shocks for (\ref{3,1}) is
described by the same self-similarity as (\ref{2.1}), with the
ODE,
 \beq
 \label{661}
  \mbox{$
 u_-(x,t)=g(z), \quad z=x/(-t)^{\frac 13} \,\,
 \Longrightarrow \,\, (g g'')'= \frac 13 \, g' z, \quad g(\mp
 \infty)= \pm 1.
  $}
   \eeq
   Existence and uniqueness for (\ref{661}) is proved similar to
   \cite[Prop.~3.1]{3NDEI}.
   In Figure
\ref{F1N}(a), we show a few  similarity profiles that create as $t
\to 0^-$ the shocks. 
 The profile for $S_-(x)$
has the derivative at the origin
 $$
 g'(0)=-0.702... \quad (\mbox{instead of $g'(0)=-0.51...$ for the NDE--3
 (\ref{1}).})
 $$
 In (b) explaining formation of $H(-x)$, the right-hand interface
 is situated at $z_0=1.297...$.
 As another known key feature,
Figure \ref{F2N} shows the  saw-type profile
  for the ODE
 \beq
 \label{663}
   \begin{matrix}
 u_-(x,t)=(-t)^\a g(z), \quad z=x/(-t)^{\frac {1+\a}3} \,\,
 \Longrightarrow \,\, (g g'')'= \frac {1+\a}3 \, g' z- \a g,
 \qquad\qquad
 \ssk\ssk\\
  \mbox{where}\quad
  \fbox{$\a_{\rm c} \approx
  -0.2384.$}\qquad\qquad\qquad\qquad
 \end{matrix}
 \eeq


The compacton equation associated with (\ref{3,1}) takes the form
 $$
 u_t=(u u_{xx})_x+ 2u u_x
 $$
 and admits the TW solution with the same $f_{\rm c}$ as in
 (\ref{62}), but now for $\l=-1$,
 \beq
 \label{789}
 u_{\rm c}(x,t)= f_{\rm c}(x+t),
 \quad
 f_{\rm c}(y)= \left\{
 \begin{matrix} 2 \cos^2( \frac y 2) \,\,\, \mbox{for} \,\,\, |y|
 \le \pi, \ssk\\
\quad\,\, 0 \qquad \,\,\, \mbox{for} \,\,\, |y|
 \ge \pi.
 \end{matrix}
  \right.
  \eeq
As for (\ref{1}), it is $\d$-entropy and G-admissible; Proposition
\ref{Pr.C} is proved similarly.

\subsection{Shock similarity profiles for Harry Dym-type
equations}

\begin{figure}
\centering
\includegraphics[scale=0.70]{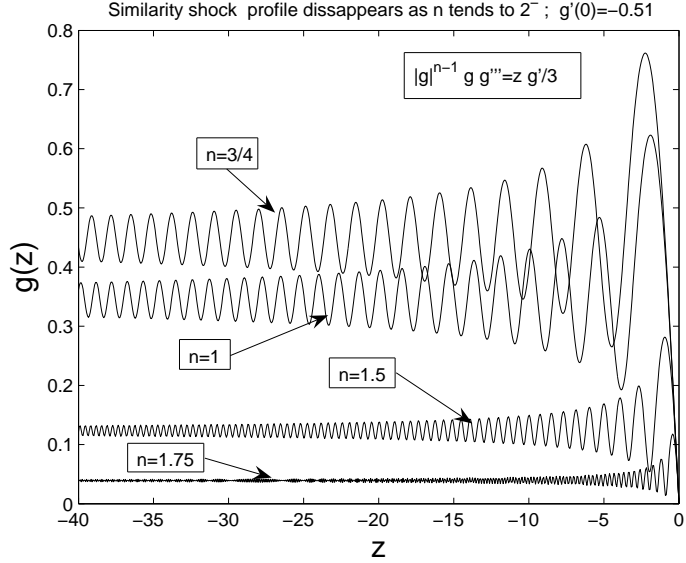}
 \vskip -.3cm
\caption{\small Shock similarity profiles of the ODE  (\ref{HD2})
for various $n \in (0,2)$; $g(0)=g''(0)=0$.}
   \vskip -.3cm
\label{FHD1}
\end{figure}

Consider the NDE
 \beq
 \label{HD1}
  u_t=|u|^{n-1}u
   \, u_{xxx} \quad (n>0),
   \eeq
   which for $n=3$ becomes the 
  quasilinear
  {\em Harry Dym}
 {\em equation}
  \beq
 \label{HD0}
  u_t = u^3 u_{xxx} \, ,
  \eeq
  which also belongs to the NDE family and
   is
an exotic integrable soliton equation; see \cite[\S~4.7]{GSVR} for
survey and references
 therein.
It admits the same formation of shocks $S_-(x)$ by the similarity
solutions given in  (\ref{661}) with the ODE
 \beq
 \label{HD2}
  \mbox{$
 |g|^{n-1}  g
 \, g'''= \frac 13 \, g' z.
 $}
  \eeq
Figure \ref{FHD1} shows that such similarity profiles exist for $n
\in (0,2)$ and vanish as $n \to 2^-$ (proof is easy), so that for
$n=3$ (the Harry Dym case) such shocks are not available.

\subsection{Shocks for fully nonlinear NDE}
 \label{SectFN}

\begin{figure}
\centering
\includegraphics[scale=0.75]{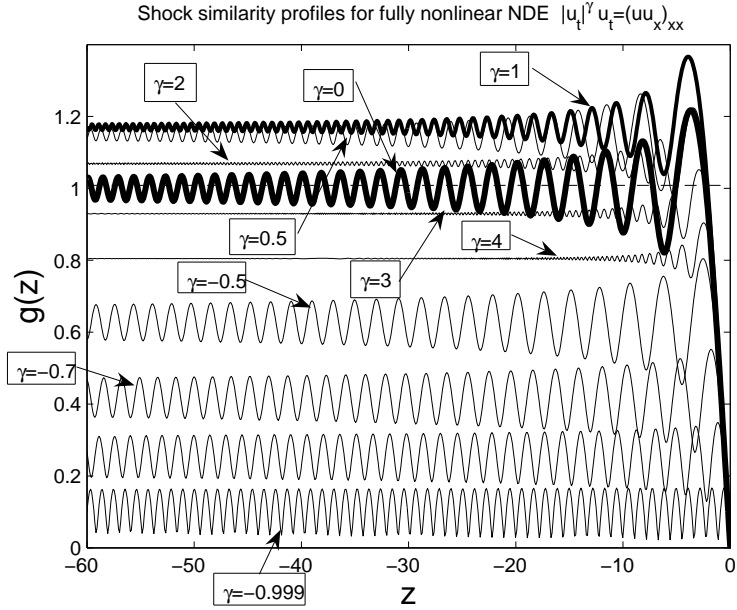}
 \vskip -.3cm
\caption{\small Shock similarity profiles of the ODE in
(\ref{FN2}) for various $\g \in (-1,4]$; $g(0)=g''(0)=0$,
$g'(0)=C=-0.51$.}
   \vskip -.3cm
\label{Ffn1}
\end{figure}

For the  NDE (\ref{FN1}), the basic blow-up similarity solutions
are slightly different,
 \beq
 \label{FN2}
  \mbox{$
 u_-(x,t)=g(z), \quad z=x/(-t)^\b, \,\,\, \b= \frac {1+\g}3
 \quad \Longrightarrow \quad (gg')''= \b^{1+\g} |g'z|^\g g'z.
 $}
  \eeq
Mathematics of such ODEs is not much different than that for
(\ref{2.1}). In Figure \ref{Ffn1}, we show how the shock
similarity profiles $g(z)$ depend on $\g > -1$. All these profiles
satisfy the anti-symmetry conditions at the origin,
 \beq
 \label{2.4}
 g(0)=g''(0)=0,
  \eeq
 and the
following expansion holds:
 \beq
 \label{FN3}
  \mbox{$
  g(z)= Cz + \b^{1+\g} \frac {C^\g}{(2+\g)(3+\g)(4+\g)}\, |z|^\g
  z^3 + ... \, .
   $}
   \eeq
   Note that the
   linearization about the constant equilibrium $C_-=1$ as $z \to
   -\infty$, again yields a nonlinear ODE,
    \beq
    \label{FN4}
    g'''= \b^{1+\g} |g'|^\g g' |z|^\g z + ... \quad (z \ll -1),
     \eeq
     which deserves further study.
  Figure \ref{Ffn1} shows that the solutions remain equally
oscillatory for all $\g>-1$, i.e., this is not a manifestation of
the oscillatory character of the linear Airy function
that occurs at a single simplest value $\g=0$ only. Thus, all ODEs
(\ref{FN2}) with $\g >-1$ contain a strong nonlinear mechanism of
oscillations about constant equilibria.

\subsection{Shock similarity profiles for cubic NDEs}
  \label{S9.3}

\begin{figure}
\centering
\includegraphics[scale=0.70]{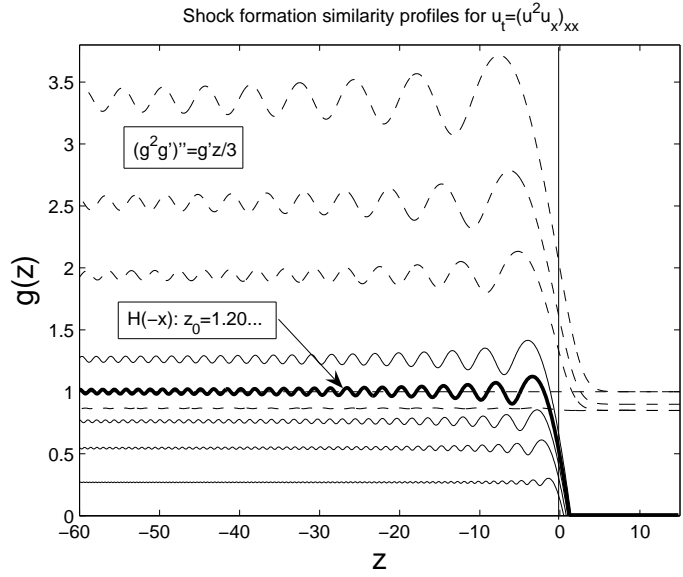}   
 \vskip -.3cm
\caption{\small Shock similarity profiles satisfying (\ref{HD2N}),
(\ref{ff1SS}) for various $z_0>0$; the boldface profile leading to
$H(-x)$ has $z_0=1.20...$; dotted lines denote shock profiles with
$g(+\infty) >0$.}
   \vskip -.3cm
\label{FHD1N}
\end{figure}

Analogously, in a similarity fashion,  the shock formation is
studies
 for the cubic fully divergent NDE (\ref{HD1N}).
The formation of shocks $H(-x)$ is described by the similarity
solutions (\ref{661}), where
 \beq
 \label{HD2N}
  \mbox{$
 (g^2 g')''= \frac 13 \, g' z,
 $}
  \eeq
  which admits a similar rigorous study.
Figure \ref{FHD1N} shows similarity profiles with the finite
interface at $z=z_0>0$ with the expansion as $z \to z_0$
 \beq
 \label{ff1SS}
  \mbox{$
  g(z)=\sqrt{\frac {z_0}6}\,(z_0-z)_+ - \frac 1{10} \frac 1{\sqrt{6
  z_0}}\, (z_0-z)_+^2+... \, ,
   $}
    \eeq
    for which the flux $(g^2 g')' \equiv \frac 13\,(g^3)''$ is continuous at $z=z_0$, so
    these are weak solutions.
    The flux is not zero for a more singular expansion such as
 $$
 \mbox{$
 g(z)= C(z_0-z)^{\frac 23}+ \frac{3z_0}{4C}(z_0-z)^{\frac43}+...
 \quad (z<z_0, \,\,\, C>0).
 $}
 $$
 Similar to \cite[\S~3.1]{3NDEI}, 
  such blow-up similarity solutions
  describe the generic formation of shock waves of the type $\sim H(-x)$ for
    (\ref{HD1N}). These solutions are entropy, which
    is proved by regular analytic approximations of the ODE as
in Section \ref{Sect6}.

 By dashed lines in Figure \ref{FHD1N}, we
denote other profiles, for which $C_+=g(+\infty)>0$, so that the
corresponding blow-up similarity solutions (\ref{2.1}) lead to
more general shocks with different values $C_\pm$ as $z \to \pm
\iy$
(with $C_+>0$). Then, as
$z \to +\infty$,
 $g(z)$ approaches $C_+$ exponentially fast,
 $$
  \mbox{$
 g(z) = C_+ + O({\mathrm e}^{-a_0 z^{3/2}}), \quad \mbox{where} \quad  a_0= \frac 2{3
 \sqrt 3 \, C_+}.
  $}
  $$

Thus,  the above solutions with the behaviour (\ref{ff1SS}) close
to interfaces show finite propagation for the NDE (\ref{HD1N}).
There are also TWs with finite
 interfaces given by
 $$
 u(x,t)=f(x +t) \quad \Longrightarrow \quad f=(f^2 f')'
 $$
that  are entropy and are approximated by the analytic family
$\{f_\d \ge \d>0\}$ satisfying
 $$
f-\d=(f^2 f')' \quad (\d>0).
 $$
 For instance, the following TW with the interface at $y=0$  is $\d$-entropy:
 $$
  f(y)=
  \left\{
 \begin{matrix}
 \frac 1{\sqrt 2}\, (-y) \,\,\, \mbox{for} \,\,\, y<0,
  \ssk\\
   \quad\,\,\, 0 \quad \,\,\, \,\, \mbox{for} \,\,\, y \ge 0.
 \end{matrix}
  \right.
 $$
 Other {\em discontinuous} TWs may not admit smooth approximations via similar
 TWs.

The boldface line in Figure \ref{FHD1N} indicates the profile that
leads to $H(-x)$ as $t \to 0^-$. Here the shock $H(-x)$ is not a
weak solution of the NDE (\ref{HD2N}). Recall that it is a
$\d$-entropy solution, i.e., there exists a converging sequence of
its smooth $\d$-deformations.

\begin{figure}
\centering
\includegraphics[scale=0.70]{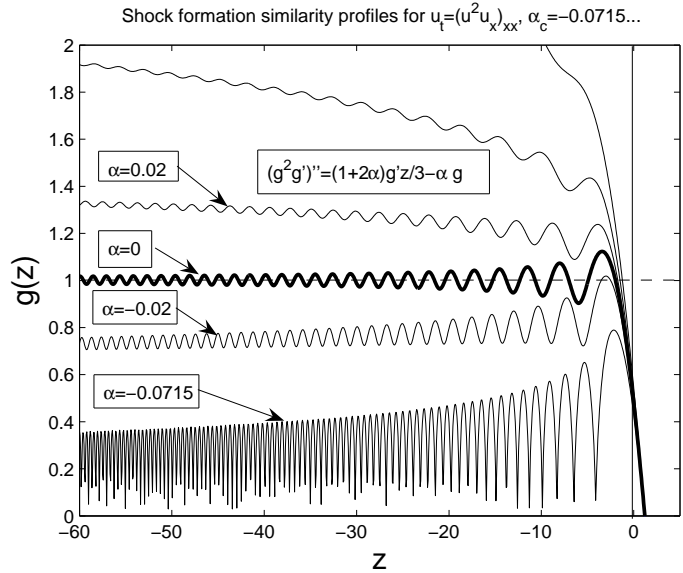}
 \vskip -.3cm
\caption{\small Shock similarity profiles of the ODE  in
(\ref{z1}) satisfying expansion (\ref{ff1N}) for $z_0=1.2$  for
various positive and negative $\a$.}
   \vskip -.3cm
\label{Fa1}
\end{figure}

More advanced shock patterns are created by  similarity solutions
(\ref{3.1}), with
 \beq
  \label{z1}
  \mbox{$
  \b= \frac {1+2 \a}3 \quad \mbox{and} \quad (g^2 g')''=
  \frac{1+2\a}3 \, g'z - \a g.
   $}
   \eeq
The interface expansion (\ref{ff1SS}) changes into
 \beq
 \label{ff1N}
  \mbox{$
  g(z)=\sqrt{\frac {(1+2\a)z_0}6}\,(z_0-z)_+ - \frac 1{10} \frac
  1{\sqrt{6(1+2\a)
  z_0}}\, (z_0-z)_+^2+... \, .
   $}
    \eeq
Figure \ref{Fa1} shows typical solutions of the ODE in (\ref{z1})
for $\a>0$ and $\a<0$. The most interesting ``saw-type" profiles
 occurs at
 $$
  \fbox{$
 \a_{\rm c} \approx -0.0715 \, .
 $}
  $$

\ssk

\noi\underline{\em On non-divergent cubic equation}.
 Consider briefly the  cubic NDE--(2,1),
 \beq
 \label{ssss1}
 u_t=(u^2 u_{xx})_x,
 \eeq
 which is similar, though it does not
admit finite propagation at the degeneracy level $\{u=0\}$. This
is seen by using TWs
 \beq
 \label{mm1}
  \mbox{$
  \begin{matrix}
 u(x,t)=f(x+t) \,\, \Longrightarrow \,\,
 f'=(f^2 f'')', \quad \mbox{so, on integration twice},
  \ssk\ssk\\
 f''= \frac 1f + \frac {C_1}{f^2}
 \,\,\Longrightarrow \,\, \frac 12(f')^2= \ln |f|- \frac{C_{1}}{f} + C_2.
 \end{matrix}
  $}
  \eeq
Setting $C_1=0$ by assuming  continuity of flux: $f^2 f''=0$ at
$f=0$, yields the ODE
 $$
 \mbox{$
 \frac 12(f')^2= \ln |f| + C_2
 $}
 $$
  that
does not allow any
 connection with the singular level
$\{f=0\}$.

\begin{figure}
\centering
\includegraphics[scale=0.70]{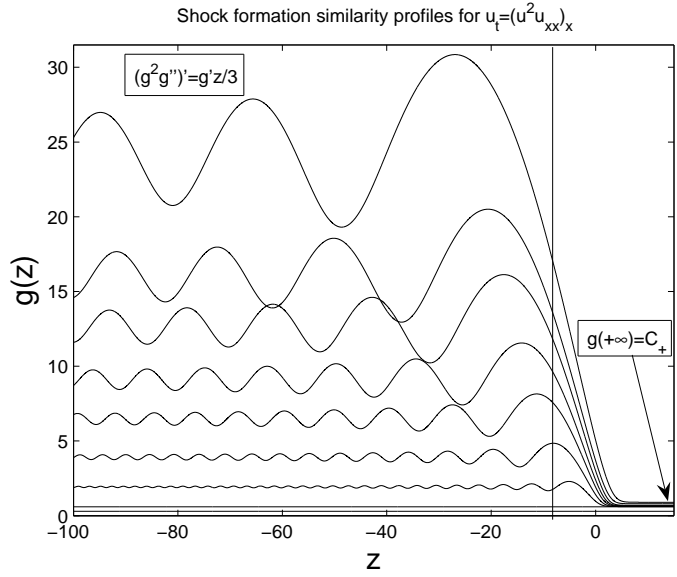}
 \vskip -.3cm
\caption{\small Shock similarity profiles satisfying the ODE
(\ref{bb111}).}
   \vskip -.3cm
\label{Fa10}
\end{figure}

The shock similarity profiles for (\ref{ssss1}) exhibit the same
form (\ref{661}) and the ODE is
 \beq
 \label{bb111}
  \mbox{$
 (g^2 g'')'= \frac 13 \, g'z.
  $}
  \eeq
Typical strictly positive profiles with
  $
  g(-\infty)=C_- > C_+=g(+\infty)>0
 $
 are shown in Figure \ref{Fa10}, so these describe blow-up
 formation of more general entropy shocks.

\begin{figure}
\centering
\subfigure[$\a=-\frac1{10}$]{
\includegraphics[scale=0.52]{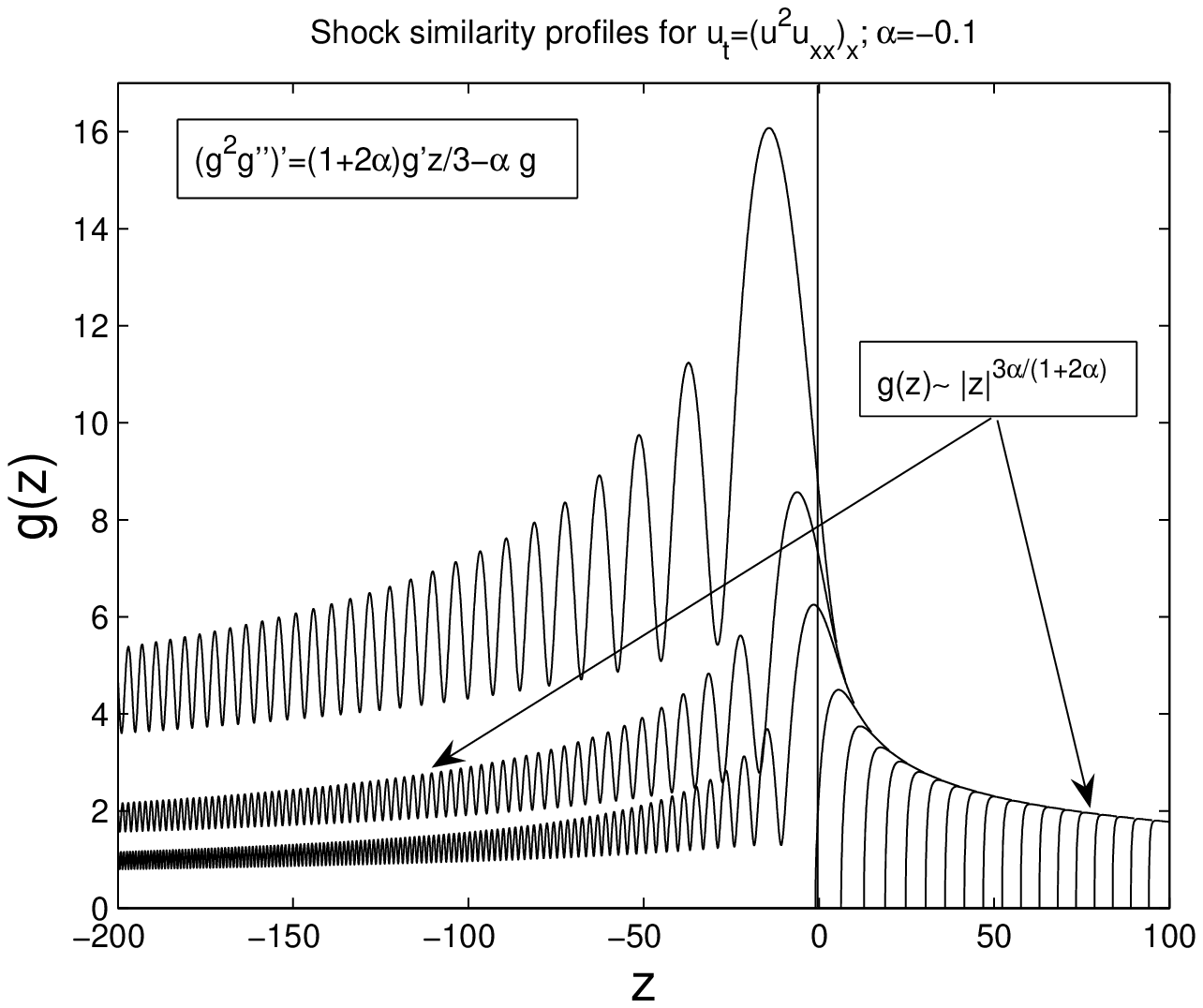} 
}
\subfigure[``Saw" for $\a_{\rm c}=-0.12559$]{
\includegraphics[scale=0.52]{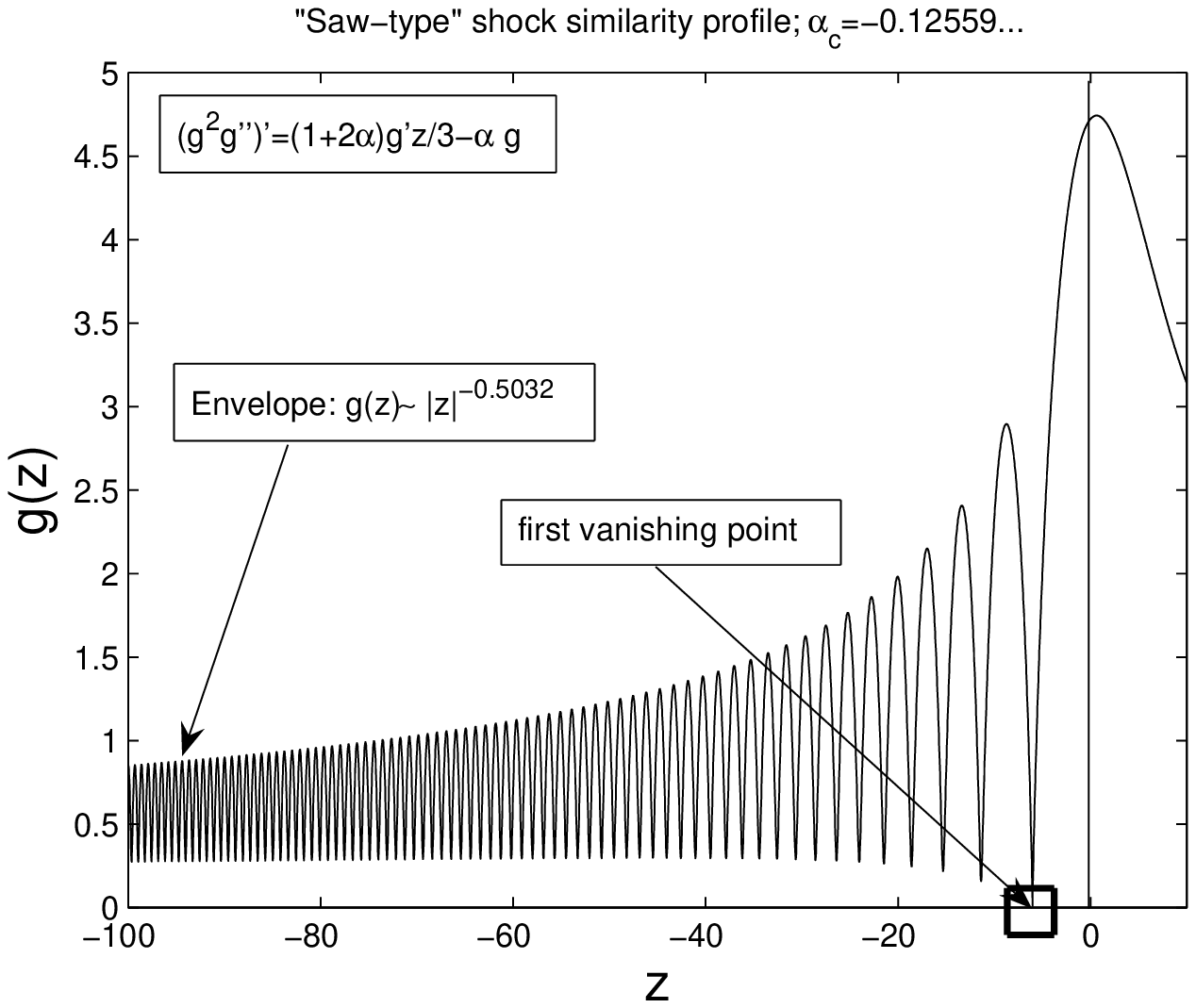} 
}
 \vskip -.4cm
\caption{\rm\small The ODE (\ref{jj1}):  the shock similarity
profiles for $\a=-\frac 1{10}$ (a), and a ``saw" profile for
$\a=\a_{\rm c}=-0.12559...$ (b).}
 \vskip -.3cm
 \label{FQQ11}
\end{figure}

More general blow-up similarity patterns (\ref{3.1}) for
(\ref{ssss1}) yields the ODE
 \beq
 \label{jj1}
  \mbox{$
 (g^2g'')'= \frac{1+2\a}3\, g'z-\a g,
  $}
  \eeq
 which exhibits properties that are similar to (\ref{z1}).
In Figure \ref{FQQ11}(a), we show typical solutions of (\ref{jj1})
for $\a=-\frac 1{10}$. These profiles are strictly positive with
 \beq
 \label{jj2}
  \mbox{$
 g(z) \sim |z|^{\frac{3\a}{1+2\a}} \quad \mbox{as}
 \quad z \to \pm \infty \quad \big(\a \in (-\frac 12,0)\big).
  $}
  \eeq
In the bottom right-hand corner of (a), we present a number of
``steep" solutions that quickly vanish (according to  (\ref{mm1})
with $C_1<0$). These show that the asymptotics (\ref{jj2}) is
unstable in the direction of shooting from $z=+\infty$.

 In (b), we present a special profile that plays a role of the
 ``saw-type" solution for
 $$
  \fbox{$
 \a=\a_{\rm c}=-0.12559... \, .
 $}
 $$
This is the best ``saw" we can get numerically, though it is seen
that there exists the first vanishing point while other ``teeth"
still stay away from zero. Anyway, we have checked that positive
shock profiles cannot be extended to $\a < \a_{\rm c}$, so this is
definitely a critical value of parameter.


\ssk

\noi\underline{\em Related compactons}. Consider the following
compacton equation
({\em q.v.} (\ref{HD1N})):
 $$
 u_t=(u^2 u_{x})_{xx}+ 9 u^2 u_x.
 $$
 The explicit compacton solution is now easier,
 $$
 u_{\rm c}(x,t)= f_{\rm c}(x+2t),
 \quad
 f_{\rm c}(y)= \left\{
 \begin{matrix} \cos y \,\,\, \mbox{for} \,\,\, |y|
 \le \frac \pi 2, \ssk\\
\quad\,\, 0 \quad  \mbox{for} \,\,\, |y|
 \ge \frac \pi 2.\,
 \end{matrix}
  \right.
  $$
Regardless the fact that it is not $C^1$ at the interface, this
solution is $\d$-entropy (note that (\ref{ff1SS}) exhibits the
same regularity). The proof uses regular approximations as in
(\ref{63}).




 \subsection{An analogy with parabolic problems}

 In a natural
sense, an analogy of the difference between the NDE--3
(NDE--(0,3)) (\ref{1}) and (\ref{3,1}) can be observed in
nonlinear parabolic theory. Namely, the fully divergent {\em
fourth-order diffusion equation}
 (the {\em DE}--4, or
{\em DE}--(0,4)),
 \beq
 \label{de1}
 u_t=- (|u|u)_{xxxx} \quad \mbox{in} \quad Q_1=(-L,L)\times (0,1)
  \eeq
     (recall that the nonlinearity $|u|u$ keeps the parabolicity on solutions of changing
sign), by classic parabolic theory \cite[Ch.~2]{LIO},  admits a
unique weak solution  of the Cauchy--Dirichlet problem with  data
$u_0$ such that $(u_0)^2 \in H^{2}$. Multiplying (\ref{de1}) by
$(|u|u)_t$ in $L^2(Q_1)$ and integrating by parts yields the
following {\em a priori} estimates of such weak solutions:
 $$
 |u|u  \in L^\infty(0,T;H^2) \quad \mbox{and} \quad (\sqrt{|u|}\, u)_t \in
 L^2(Q_1).
 $$
 Uniqueness follows from the {\em monotonicity}
of the operator $ - (|u|u)_{xxxx}$ in $H^{-2}$: for two weak
solutions $u$ and $v$,
 \beq
 \label{de2}
\begin{matrix}
 \frac 12\,  \frac{\mathrm d}{{\mathrm d}t}\|u-v\|_{H^{-2}}^2=
- \int (|u|u - |v|v)_{xxxx}\, (D_x^4)^{-1}(u-v)\qquad\ssk\ssk\ssk
\\
     = \, -
 \int(|u|u- |v|v)(u-v) \le 0,
 \end{matrix}
 \eeq
so that
  (\ref{de2})
  guarantees {continuous dependence
   of solutions
on initial data}.


On the other hand, the fourth-order {\em thin film equation}
(TFE--4)
 \beq
  \label{de3}
 u_t= - (|u| u_{xxx})_x,
  \eeq
  which has the distribution of the derivatives (3,1),
  does not admit such a simple treatment of continuous dependence
  and uniqueness as via (\ref{de2}). The Cauchy problem for the
  non-fully divergent TFE--4 (\ref{de3}) needs special
  approximation approaches, \cite{Gl4}.



 For non-fully divergent operators such as in
(\ref{3,1}) or fifth-order ones of the types (2,3), (3,2), (4,1),
in the NDEs (see \cite{GalNDE5})
 $$
u_t=-(u u_{xx})_{xxx}, \quad u_t=-(u u_{xxx})_{xx}, \quad u_t=-(u
u_{xxxx})_{x},
 $$
we face a difficulty that is similar to that for the TFE
(\ref{de3}). In both cases, the $\d$-approximation concepts
will play a role, quite similarly
to the higher-order parabolic TFEs--6 such as (see \cite{GBl6} and
references therein)
 $$
u_t=(|u| u_{xxxxx})_x, \quad u_t=(|u| u_{xxxx})_{xx}, \quad
u_t=(|u| u_{xxx})_{xxx}, \quad \mbox{etc.}
   $$


\section{On related higher-order in time NDEs}
 \label{Sect2t}

It is principal for PDE theory to justify that the ideas
       of similarity  shock
wave formation remain valid for other NDEs that are higher-order
in time. We  claim that the concept of smooth $\d$-deformations
 can be developed for such quasilinear degenerate PDEs.
 Let us present a few comments in these directions.

\subsection{Second-order in time NDE}

\begin{figure}
\centering
\includegraphics[scale=0.70]{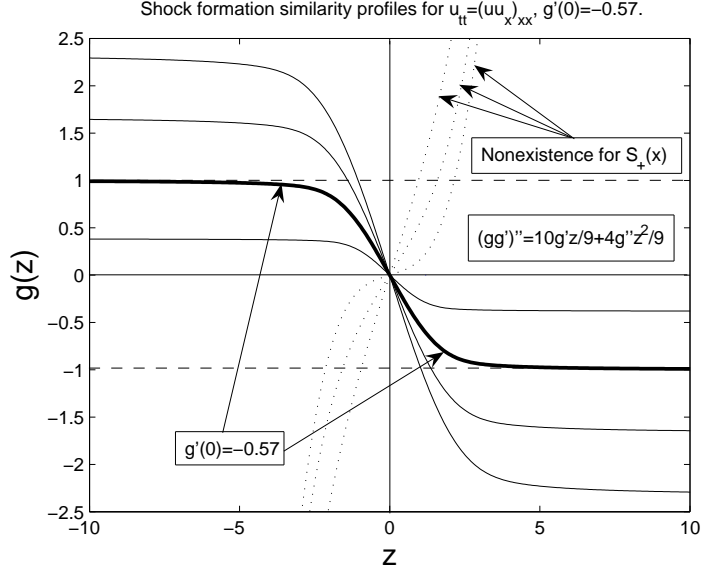}
 \vskip -.3cm
\caption{\small The shock similarity profile as the unique
solution of the problem (\ref{2.2t}).}
\label{F1t}
\end{figure}

\begin{figure}
\centering
\subfigure[various $z_0>0$]{
\includegraphics[scale=0.52]{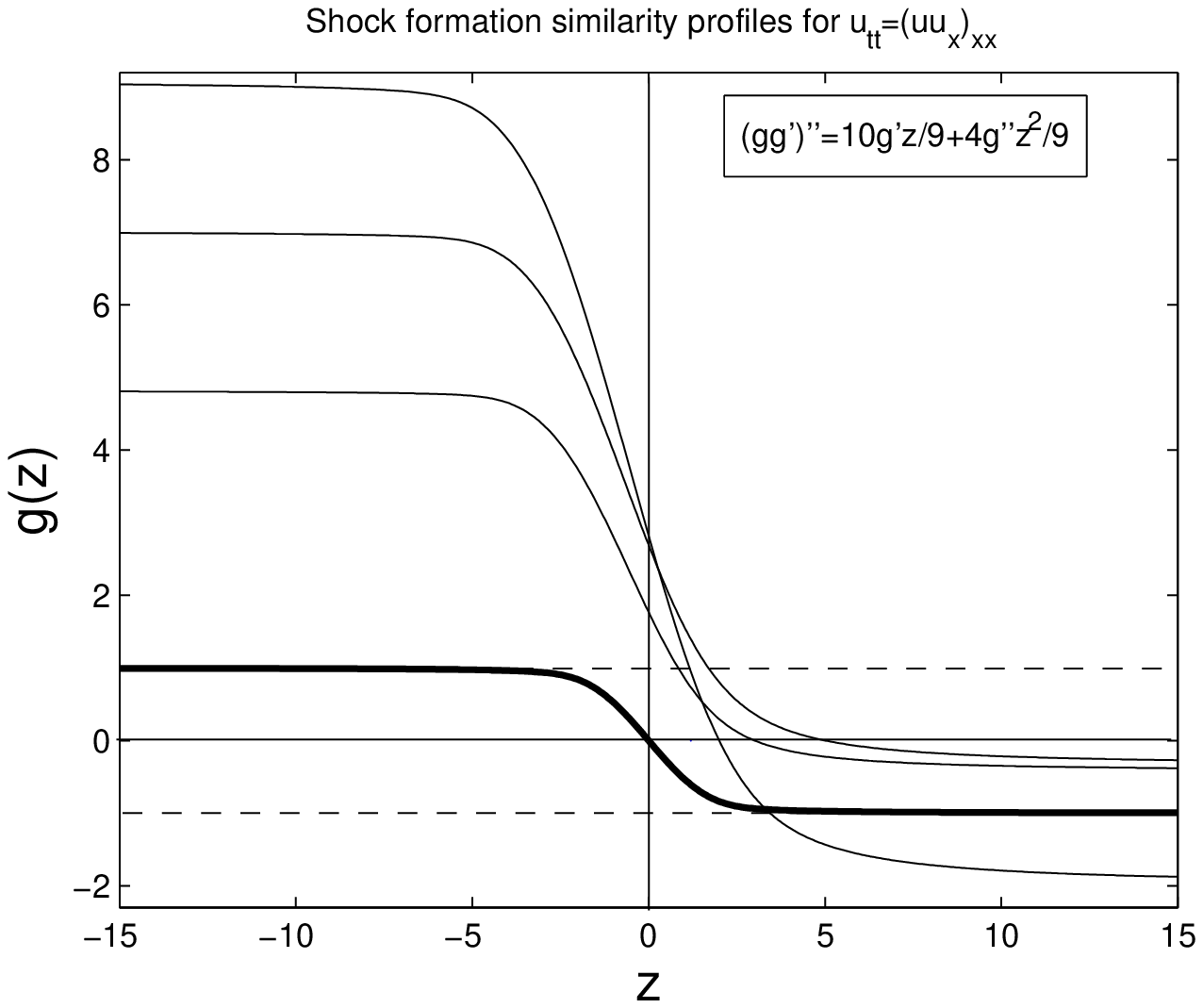} 
}
\subfigure[$z_0=5$]{
\includegraphics[scale=0.52]{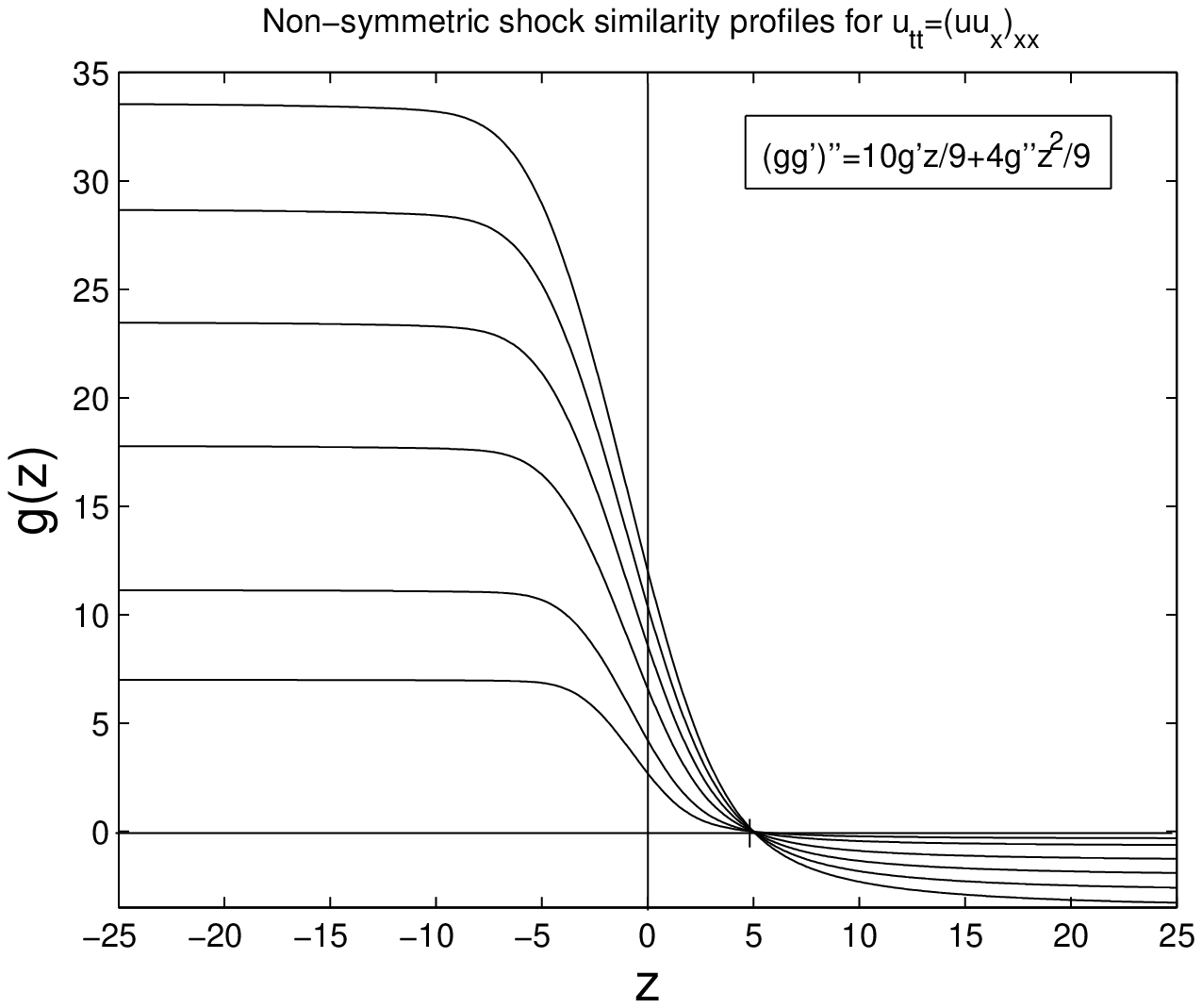} 
}
 \vskip -.4cm
\caption{\rm\small Non-symmetric shock profiles  satisfying the
ODE in (\ref{2.2t}).}
 \vskip -.3cm
 \label{F2t}
\end{figure}

 As in \cite[\S~1.2]{3NDEI},
  we begin with the simple observation:
  $S_\pm(x)$ are stationary weak solutions of the {\em
 second-order in time NDE}
  \beq
  \label{t1}
  u_{tt}=(u u_x)_{xx}.
   \eeq
   To distinguish the entropy one, as usual, we introduce  the  similarity
   solutions
 \beq
 \label{2.1t}
 u_-(x,t)=g(z), \quad z= x/(-t)^{\frac 23}, \quad \mbox{where}
  \eeq
 \beq
 \label{2.2t}
  \mbox{$
  (g g')''= \frac {10}9 \, g'z+ \frac 4 9 \, g'' z^2 \quad \mbox{in} \quad \re, \quad f(\mp
  \infty)=\pm 1.
   $}
   \eeq
The study of this ODE is similar to that in \cite[\S~3]{3NDEI},
so we present the existence result for the shock $S_-(x)$ in
Figure \ref{F1t}. The dotted lines show nonexistence of similarity
profiles for $S_+(x)$ (cf. a proof below). The boldface profile is
unique and satisfies the anti-symmetry conditions at the origin
(\ref{2.4}). We see that profiles $g(z)$ are now non-oscillatory
about $\pm 1$ and the convergence to these constant equilibria is
exponentially fast,
 $$
 g(z) = \pm 1 + O({\mathrm e}^{-\frac 4{27}\, |z|^3}) \quad
 \mbox{as}
 \quad z \to \mp \infty.
 $$
This reflects the fact that the fundamental solutions of the
corresponding  linear PDE
 \beq
 \label{ggg1}
 u_{tt}=u_{xxx}
 \eeq
 is not oscillatory as $x \to \infty$.
 Obviously, the blow-up similarity solution (\ref{2.1t}) generates
 in the limit $t \to 0^-$ the shock $S_-(x)$, i.e.,
 \beq
 \label{2.3}
  u_-(x,t) \to S_-(x) \quad \mbox{as}
  \quad t \to 0^-.
   \eeq
 In Figure \ref{F2t}, we show various non-symmetric shock
similarity profiles with different limits
 as $z \to \pm \iy$.

Incidentally,  $S_+(x)$ cannot be obtained in such a limit, since
the ODE (\ref{2.2t}) does not admit suitable similarity profiles
$g$. This can be seen from the identity obtained by multiplying
the ODE (\ref{2.2t}) by $g'$ and integrating over $(0,\infty)$
with conditions (\ref{2.4}),
 $$
  \mbox{$
  -\frac 23 \, (g'(0))^3- \int\limits_0^\infty g (g'')^2 = \frac 23 \, \int\limits_0^{\infty}(g')^2 z
  >0.
  $}
  $$
Therefore, for  $g'(0)>0$, there is no positive solution $g(z) \to
+1$ as $z \to +\infty$, since the left-hand side is then strictly
negative.

As an important difference with the previously studied PDEs,  the
NDE (\ref{t1}) is symmetric under the time-reflection,
  \beq
  \label{time1}
 t \mapsto -t
 \eeq
so that the stationary shock $S_-(x)$ that appears as $t \to 0^-$
according to the similarity law (\ref{2.1t}) (according to
``centre/stable manifold" behaviour in \cite[\S~6]{3NDEI})
     will
next disappear in the same smooth similarity manner (\ref{2.1t}),
where $(-t)$ is replaced by $t$. As a next step, the concept of
smooth $\d$-deformations should be applied to (\ref{t1}) to
produce a unique solution of the Cauchy problem, but this demands
extra more technical study.

\subsection{Third-order in time NDE}

Consider the {\em
 third-order in time NDE}
  \beq
  \label{t1tt}
  u_{ttt}=(u u_x)_{xx}, \quad \mbox{or}
   \eeq
 \beq
 \label{jj1gg}
 \left\{
 \begin{matrix}
 u_t=v_x, \,\,\,\,\\
 v_t=w_x, \,\,\,\,\\
 w_t=uu_x,
  \end{matrix}
  \right.
   \eeq
   being
 a first-order system with the characteristic
equation $\l^3=u$, with one real and two complex eigenvalues for
$u \not = 0$, so it not hyperbolic.

Quite analogously,
  $S_\pm(x)$ are stationary weak solutions of (\ref{t1tt})
  for which the basic (with $\a=0$)   similarity
   solutions are
 \beq
 \label{2.1tt}
 u_-(x,t)=g(z), \quad z= x/(-t), \quad \mbox{where}
  \eeq
 \beq
 \label{2.2tt}
  \mbox{$
  (g g')''=(z^3g')'' \equiv 6 g'z+6 g'' z^2+ g''' z^3 \quad \mbox{in} \quad \re, \quad f(\mp
  \infty)=\pm 1.
   $}
   \eeq
Integrating (\ref{2.2tt}) twice yields
 $$
 gg'=z^3 g' +Az+B, \quad \mbox{with constants} \quad A, \, B \in \re,
 $$
 so that the necessary similarity profile $g(z)$ solves
 the first-order ODE
 \beq
 \label{ss1GG}
 \mbox{$
  \frac{{\mathrm d}g}{{\mathrm d}z}= \frac {A z}{g-z^3}, \quad
  \mbox{where}
  \quad A=(g'(0))^2>0.
  $}
  \eeq
By the phase-plane analysis of (\ref{ss1GG}), we easily get the
following:

\begin{proposition}
 \label{Pr.3}
 The problem $(\ref{2.2tt})$ admits a unique solution $g(z)$
 satisfying the anti-symmetry conditions $(\ref{2.4})$ that is
 positive  for $z<0$, monotone decreasing, and is real analytic.
  \end{proposition}

\begin{figure}
\centering
\includegraphics[scale=0.70]{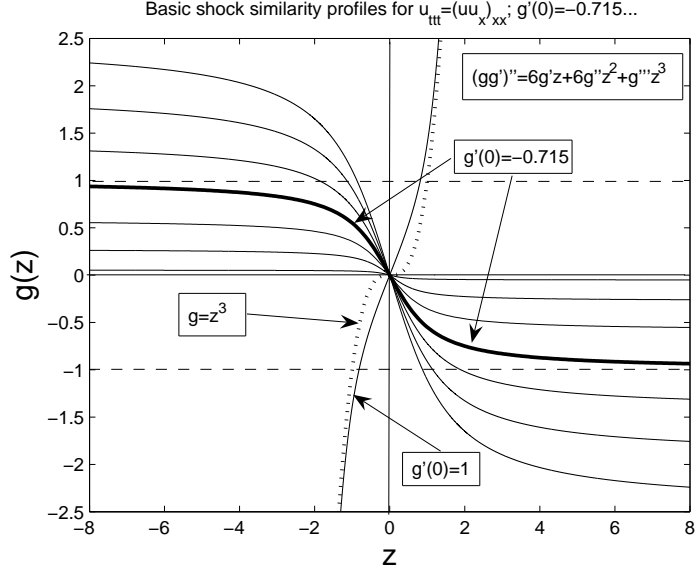}
 \vskip -.3cm
\caption{\small The shock similarity profile as the unique
solution of the problem (\ref{2.2tt}).}
   \vskip -.3cm
\label{F1tt}
\end{figure}

Such basic anti-symmetric similarity profiles  are shown in Figure
\ref{F1tt}. These satisfy  the expansion near the origin, as $z
\to 0$,
 \beq
 \label{rr1t}
 \mbox{$
 g(z)= \sum\limits_{(k \ge 0)}c_k z^{2k+1}= C z + \frac 1{4} \, z^3 + \frac{3}{32}\, \frac 1 C \,
 z^5+...  \quad (C =g'(0)<0).
  $}
  \eeq
Substituting the expansion in (\ref{rr1t}) into (\ref{ss1GG})
yields
 $$
  \mbox{$
  g'(g-z^3)=Az \,\,\, \Longrightarrow \,\,\,
  \sum\limits_{(k,j\ge 0)} (2k+1)c_k(c_j-\d_{j1}) z^{2(k+j)+1}=Az.
   $}
   $$
The corresponding algebraic system for the expansion coefficients
$\{c_k\}$ is uniquely solved giving the unique analytic solution.
The boldface profile $g(z)$ in Figure \ref{F1tt} (by (\ref{2.1tt})
it gives $S_-(x)$ as $t \to 0^-$) is non-oscillatory about $\pm 1$
with the algebraic convergence
 $$
  \mbox{$
 g(z) = \pm 1 + \frac A{z} +... \quad
 \mbox{as}
 \quad z \to \mp \infty.
  $}
 $$
Again, the fundamental solutions of the corresponding linear PDE
 \beq
 \label{dd1}
 u_{ttt}=u_{xxx}
 \eeq
 is not oscillatory as $x \to \pm \infty$.
  The linear PDE
  (\ref{dd1}) exhibits  some finite propagation features with the corresponding
 test consisting of  checking the TWs,
 $$
 u(x,t)=f(x- \l t) \quad \Longrightarrow \quad
  -\l^3 f'''=f''', \,\,\,\, \mbox{i.e.,} \,\,\, \l=-1,
   $$
   where the profile $f(y)$ disappears from. This is similar to
   a few other well-known canonical equations of
   mathematical physics such as
    $$
    u_t=u_x \,\,(\mbox{dispersion, $\l=-1$}) \quad \mbox{and} \quad
  u_{tt}=u_{xx} \,\,(\mbox{wave equation, $\l=\pm 1$}).
   $$
 Any finite propagation is not true for (\ref{ggg1}).
 The blow-up  solution (\ref{2.1tt}) gives
 in the limit $t \to 0^-$ the shock $S_-(x)$, and (\ref{2.3})
 holds.  In Figure \ref{F1tt} we also show the results of shooting
 with $g'(0)>0$ giving unbounded profiles $g(z) \sim z^3$ as $z
 \to \pm \infty$. As usual, this means nonexistence of similarity blow-up
 profiles corresponding to $S_+$-type shocks.

A key difference with the previous  problems is that the original
ODE (\ref{2.2tt}) written as
 \beq
  \label{lc}
\mbox{$
  (g-z^3)g'''= 6 g'z+6 g'' z^2-3g'g''
   $}
   \eeq
has, instead of $\{g=0\}$, another {\em singular line} (a kind of
nonlinear ``{\em light cone}")
 \beq
 \label{lc1}
 L_0: \quad g(z)=z^3.
 \eeq
 Then, formally, the existence of global solutions of (\ref{lc}) depends on the
 possibility of a continuous transition through it.
 The simpler  integrated form (\ref{ss1GG}) shows  that
 typical solutions do not cross $L_0$ (except at the analytic point $z=0$), so that
  ``weak discontinuities" do no occur.

Since (\ref{t1tt}) has the same symmetry
 \beq
 \label{symm1}
 \left\{
 \begin{matrix}
u \mapsto -u, \\
 t \mapsto -t, \,
  \end{matrix}
  \right.
  \eeq
  as
(\ref{1}), similarity solutions (\ref{2.1tt}) with $-t \mapsto t$
and $g(z) \mapsto g(-z)$
  also give the
rarefaction waves for $S_+(x)$, as well as other types of collapse
of initial non-entropy discontinuities.



Using the known asymptotic properties of blow-up similarity
solutions (\ref{2.1tt}) and those global with $-t \mapsto t$, for
convenience, we formulate the following

\begin{proposition}
 \label{Pr.CK1}
 The Cauchy problem for
 the equation $(\ref{t1tt})$  admits:

 \noi{\rm (i)} an analytic solution $u_-(x,t)$ in $Q_T=\re \times (0,T)$ that
 converges as $t \to T^-$ to the shock $S_-(x)$ in $L^1_{\rm loc}$
 and  a.e., and

\noi{\rm (ii)} for non-analytic singular initial data as $t \to
0^+$ given by
 \beq
 \label{data1}
  \mbox{$
 u(x,t) \to S_+(x), \quad u_t(x,t) \to \frac Ax, \quad \frac 1{t^2}\,
 u_{tt}(x,t)\to
 \frac {3A \,{\rm sign} \, x}{x^4}
  $}
   \eeq
with uniform convergence as $t \to 0$ on any compact subset from
$\re \setminus\{0\}$ $($and in $L^1_{\rm loc}$ for $u(x,0))$,
   there exists an analytic
solution in $\re \times \re_+$.

\end{proposition}

\smallskip

\noi\underline{\em Analytic $\d$-deformations by
Cauchy-Kovalevskaya theorem}.
Eventually, we start to deal with the third-order in time NDE
(\ref{t1tt}) that turns out to be in the {\em normal form}, so it
obeys the Cauchy--Kovalevskaya (C-K) theorem \cite[p.~387]{Tay}.
Hence, for any analytic initial data $u(x,0)$, $u_t(x,0)$, and
$u_{tt}(x,0)$, there exists a unique local in time analytic
solution $u(x,t)$\footnote{In this connection, the result (ii) in
Proposition \ref{Pr.CK1} sounds unusual: for non-analytic and very
singular data, there exists a global analytic solution.}. Thus,
(\ref{t1tt}) generates a local semigroup of
 analytic solutions, and this makes it easier to deal with smooth
 $\d$-deformations that always can be chosen to be analytic.
  On the other hand, such nonlinear PDEs can
 admit other (say, weak) solutions that are not analytic.
Actually, Proposition \ref{Pr.3} shows that the shock $S_-(x)$
 is a $\d$-entropy solution of (\ref{t1tt}), which is obtained by finite-time blow-up
  as $t \to 0^-$
 from the analytic similarity solution (\ref{2.1tt}).

\ssk

\noi\underline{\em Shocks  for  non-degenerate NDE}. For the
corresponding non-degenerate NDE
 \beq
 \label{lk1}
 u_{ttt}= ((1+u^2)u_x)_{xx},
  \eeq
  the similarity solutions (\ref{2.1tt}) lead, on double
  integration, to the ODE (cf. (\ref{ss1GG}))
   \beq
   \label{lk2}
   \mbox{$
   g'= \frac{Az+B}{1+g^2-z^3} \quad (A,\, B \in \re).
    $}
    \eeq
It is easy to show using the phase-plane, that for $z_0= - \frac
BA
>1$ (this gives a necessary extra singular point $(z_0,g_0)$ of the flow, where
 $g_0^2= z_0^3-1$), (\ref{lk2}) admits  analytic solutions $g(z)$
satisfying $g(\pm \infty)=C_\pm>0$ with $C_->C_+$, so as $t \to
0^-$, we obtain the shock.

\subsection{Stationary entropy shocks for other higher-order in time NDEs}

We now very briefly check entropy properties of the shocks
$S_\pm(x)$ for the following NDEs of arbitrary order:
 \beq
 \label{aa1}
 D_t^{2m+2} u=D_x^{2m}(u u_x) \quad (m \ge 1).
  \eeq
For $m=0$, this gives the following simple NDE:
 \beq
 \label{aa2}
 u_{tt}=u u_x,
  \eeq
  for which both shocks $S_\pm$ are obviously weak solutions, so
  one needs to identify which ones are entropy.
Note that, as (\ref{t1tt}), the PDEs (\ref{aa1}) for any $m \ge 1$
obey the Cauchy--Kovalevskaya theorem, so a unique local semigroup
of analytic solutions does exist.

\ssk


\noi\underline{\em $\d$-entropy $S_-$ via analytic TWs}.
 For a change, we  present $\d$-deformations by TWs
  \beq
  \label{aa3}
   \begin{matrix}
   \mbox{$
  u(x,t)=f_\l(x-\l t) \quad \Longrightarrow \quad
 \l^{2m+2} f^{(2m+2)}= \frac 12\,(f^2)^{(2m+1)}, \quad \mbox{or}
  $}
   \ssk\\
  \mbox{$
\l^{2m+2}f'=-\frac 12\,(1-f^2) \quad \Longrightarrow \quad
   f_\l(y)=\frac{{\mathrm e}^{-y/\l^{2m+2}}
-1}{{\mathrm e}^{-y/\l^{2m+2}}+1}.
 $}
  \end{matrix}
 \eeq
We then observe that
 \beq
 \label{aa5}
  f_\l(y) \to S_-(y) \quad \mbox{as \,\, $\l \to 0$\,\, uniformly
  in $\re$,}
   \eeq
   so that the stationary shock wave $S_-(x)$ is G-admissible and
   is $\d$-entropy, where the necessary $\d$-deformation is given
   by the TW (\ref{aa3}) with $\l=\d$.

A similar (but not explicit) construction of $\d$-entropy
solutions with convergence (\ref{aa5}) is performed for other
normal NDEs such as
 \beq
 \label{aa10}
 D_t^{2m+4} u=-D_x^{2m}(u u_x)
, \quad \mbox{or} \quad  D_t^{2m+2k} u=(-1)^{k+1}D_x^{2m}(u u_x),
\,\,\, k \ge 1.
  \eeq
The corresponding analytic TW profiles  $f_\l(y)$ satisfying the
convergence  (\ref{aa5}) in $L^1_{\rm loc}$ are described in
\cite[\S~4]{GalEng}.

\ssk

\noi{\bf Remark: $S_+$ can be formally created by a classical but
non-analytic blow-up self-similar solution}. There exists a
self-similar blow-up to $S_+$ for the NDE (\ref{aa2}) via
 \beq
 \label{aa6}
  \mbox{$
 u_-(x,t)=g(z), \,\, z= \frac x{(-t)^2} \,\, \Longrightarrow \,\, g'z^2=
 \frac 14 \,g^2, \,\, \mbox{so}
  \,\,\,
   g(z)= \left\{
  \begin{matrix}
  \frac{4z}{4z+1}, \,\, z \ge 0, \ssk\\
 \frac{4z}{1-4z}, \,\, z \le 0.
  \end{matrix}
  \right.
  $}
  \eeq
  This ODE obeys the symmetry
 \beq
     \label{symm88}
      \left\{
      \begin{matrix}
     g \mapsto -g,\\
     z \mapsto -z,
      \end{matrix}
      \right.
      \eeq
  Note that $g(z)$ is just $C^1$ (not $C^2$) at $z=0$, which is
   enough to represent a weak solution of the degenerate PDE (\ref{aa2})
   (though, as we know, being weak often means almost nothing).
   Moreover, (\ref{aa6}) is a classical $C^{1,2}_{x,t}$ solution
   of (\ref{aa2}).
    Observe that the non-analyticity of $g(z)$ is
    associated with the too strong degeneracy at $z=0$ of the
   corresponding ordinary differential operator
    $
    \mbox{$
    z^2 \, \frac{\mathrm d}{{\mathrm d}z}.
     $}
     $
We  suspect that (\ref{aa6}) is not entropy at all. Moreover, one
can see that $g(z)$ {\em is not} an odd function, so it looks more
like a solution of an IBVP for $x>0$ with some boundary condition
at $x=0$. Nevertheless, we recall that it is a classical solution
of the Cauchy problem.
 The NDEs (\ref{aa1}) deserve deeper study.




\section{On shocks for spatially higher-order  NDEs}
 \label{Sect55}

\subsection{Fifth-order NDEs}

\begin{figure}
\centering
\includegraphics[scale=0.75]{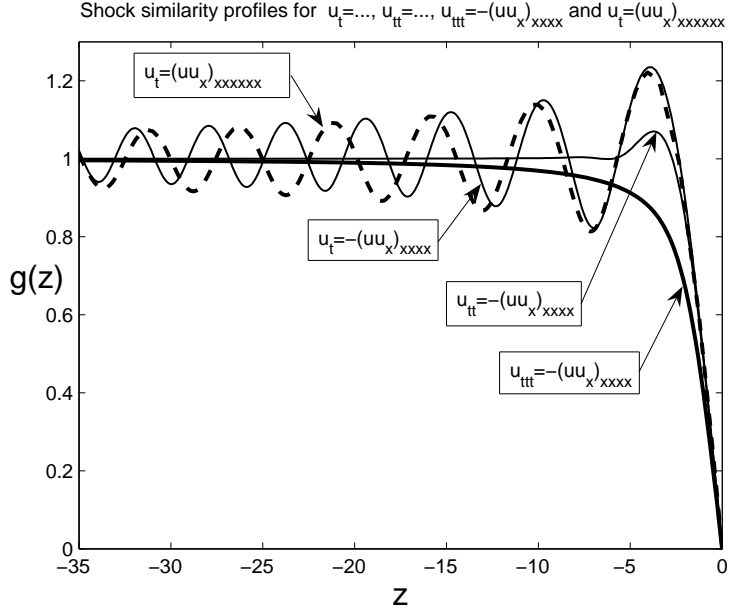}
 \vskip -.3cm
\caption{\small The shock similarity profiles $g(z)$ for $z<0$
 as
the unique solutions of the problems (\ref{55}), (\ref{56})  and
(\ref{771}).}
   \vskip -.3cm
\label{F55t}
\end{figure}

The similarity mechanism of shock formation remains valid for
higher-order NDEs, among which, as an illustration, we comment on
the following three (including the NDE--5 (\ref{N5})):
 \beq
 \label{51}
 \begin{matrix}
 u_t=-(uu_x)_{xxxx},\ssk\,\,\\
 u_{tt}=-(uu_x)_{xxxx},\ssk\,\\
 u_{ttt}=-(uu_x)_{xxxx}.
 \end{matrix}
 \eeq
 Concerning application of such fifth and higher-order NDEs, see
 \cite{GalNDE5},
 \cite[p.~166]{GSVR}, and references therein.
The blow-up similarity solutions of $S_-$-type are the same,
 \beq
 \label{54}
 \mbox{$
 u_-(x,t)=g(z), \quad z=x/(-t)^\b, \quad \mbox{where} \quad \b=
 \frac 15, \quad \b= \frac 25, \quad \b= \frac 35,
  $}
  \eeq
  respectively. The ODEs are, respectively,
   \beq
   \label{55}
   \begin{matrix}
   (gg')^{(4)}= - \frac 15 \, g'z,\qquad\qquad\qquad\qquad\qquad\ssk\\
   (gg')^{(4)}= - \frac 2{25}(7 g'z+2 g'' z^2),\qquad\qquad\quad\ssk\\
    (gg')^{(4)}= - \frac 3{125}(64 g'z+57 g'' z^2+9 g''' z^3).
    \end{matrix}
    \eeq
 These are much more complicated equations than all those studied before.
 We  do not have a proof of existence of the $S_-$-type
 profiles $g(z)$ to say nothing about uniqueness, though we can justify that the shooting
 procedure to get a solution is well-posed according to dimensions of stable and unstable manifolds of
 orbits at the singular points $z=0$ (where $g=0$) and $z=-\infty$ (where $g=+1$). On the other
 hand, the same numerical methods give us a strong evidence of
  existence-uniqueness. In Figure \ref{F55t}, using  {\tt bvp4c} solver of {\tt MatLab}, we
 present the unique solutions of the ODEs (\ref{55}) satisfying
 the standard conditions
  \beq
  \label{56}
  g(\pm \infty)=\mp 1 \quad \mbox{and} \quad
  g(0)=g''(0)=g^{(4)}(0)=0 \,\,\,(\mbox{anti-symmetry}).
   \eeq
   Note that first two ODEs admit solutions that are oscillatory
   about the equilibrium $g=1$ as $z \to - \infty$, while the last
   one has monotone non-oscillatory solutions according to the
   following asymptotics, respectively: for $z \ll -1$, neglecting
   lower-order algebraic multipliers in the second and third
   formulae,
   \beq
   \label{57}
   \begin{matrix}
g(z)-1 \sim |z|^{-\frac 38} \cos\big(\frac 4{5 \sqrt 2}\,5^{-\frac
14}|z|^{\frac 54}+c_0\big), \qquad\,\,\,\,\,\,\,\ssk\\
 g(z)-1 \sim {\mathrm e}^{-\frac 3{5\sqrt 2}(\frac
 25)^{\frac 23}|z|^{\frac 53}} \cos\big( \frac 3{5\sqrt 2}(\frac
 25)^{\frac 23}|z|^{\frac 53}+c_0\big), \ssk\\
 g(z)-1 \sim {\mathrm e}^{- \frac 25(\frac 35)^{\frac
 32}|z|^{\frac 52}}.\qquad\qquad\qquad\qquad\,\,\qquad
  \end{matrix}
   \eeq
   The exponentially small oscillations in the second
   line are hardly seen in the figure and requires another,
   logarithmic scale for revealing those.

\subsection{On a seventh-order NDE}

For completeness and convenience of comparison,  Figure \ref{F55t}
also gives the shock similarity profiles (the dashed line) for the
NDE--7,
 \beq
 \label{7eq}
 u_t=(u u_x)_{xxxxxx}, \quad \mbox{where}
 \eeq
  \beq
  \label{771}
  \mbox{$
 u_-(x,t)=g(z), \quad z=x/(-t)^{\frac 17}
 \quad \Longrightarrow \quad (g g')^{(6)}= \frac 17 \, g'z, \quad
 g(\pm \infty)=\mp 1.
 $}
 \eeq
The shock profile is very similar to that for the NDE--5 in
(\ref{N5}), so that a general geometry of these  shock profiles
does not essentially depend on the order, $(2m+1)$, of the PDEs
(\ref{N5}) for $m \ge 1$; the oscillatory behaviour also changes
slightly with $m$ and always has the type given in the first line
in (\ref{57}).

These results show that, for all the above higher-order NDEs,
canonical shocks of $S_-$-type are obtained by blow-up in finite
time from smooth classical solutions. According to our
$\d$-entropy approach, this confirms a correct entropy nature of
such shock waves.

\begin{figure}
\centering
\includegraphics[scale=0.85]{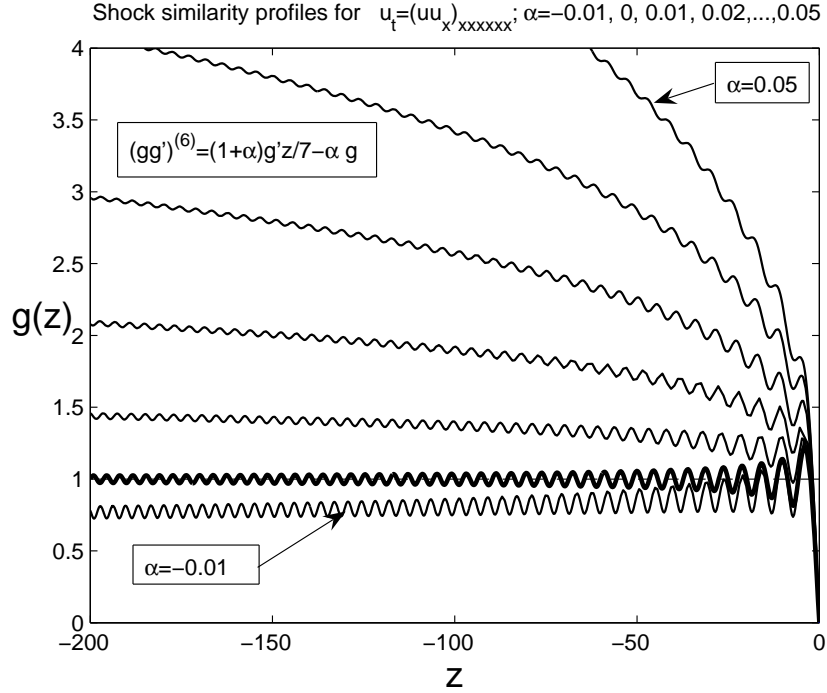}
 \vskip -.3cm
\caption{\small Similarity profiles $g(z)$ for $z<0$
 satisfying the ODE in (\ref{771N}) for various $\a \in [-0.01,0.05]$.}
   \vskip -.3cm
\label{F55tN}
\end{figure}

Let us describe other types of shocks and rarefaction waves for
(\ref{7eq}) driven by blow-up similarity patterns
 \beq
  \label{771N}
  \mbox{$
 u_-(x,t)=(-t)^\a g(z), \quad z=x/(-t)^{\b}, \,\,\, \b=
 \frac{1+\a}7
 \quad \Longrightarrow \quad (g g')^{(6)}= \frac {1+\a}7 \, g'z- \a g.
 $}
 \eeq
 These similarity profiles are presented in Figure \ref{F55tN}.
This shows that the profiles get more oscillatory for
$\a<0$, but we failed to detect a ``saw"-type profile as in
 \cite[\S~4.3]{3NDEI}
for such a seventh-order ODE by using any numerical method.

Finally,  the analysis of the ODE in (\ref{771N}) on the invariant
subspace (cf. the invariant subspace in \cite[\S~4.3]{3NDEI})
 $
 W_4={\rm Span}\{z,z^3,z^5,z^7\}
 $
 shows that a nontrivial dynamics  exists for the critical exponent
  $$
   \fbox{$
  \mbox{$
  \a_{\rm c}= \frac{415}{2574}=0.161228... \, ,
   $}
    $}
   $$
   and that the explicit solutions are given by
   $$
    \mbox{$
   g(z)= C z + \frac{6!}{13!}\, z^7, \quad \mbox{where}
   \quad C \in \re\,\, \mbox{is arbitrary}.
 $}
    $$

\section{On changing sign compactons for higher-order NDEs}
 \label{SectComp5}

Finally, we return to the compacton solutions of the NDEs.
First time, we discussed
 the entropy properties of compactons in Section \ref{Sect6} for
 the NDE--3, where the entropy nature of such solutions was
 successfully justified.
 It turns out that the fact that these compactons are $\d$-entropy,
 i.e., are constructed by smooth $\d$-deformations,
 can be proved by a purely ODE approach, by smooth positive
 approximations of compactons via analytic solutions.
We must admit that this ODE approach cannot be extended in
principle to higher-order NDEs, so we need either to return to the
original PDE $\d$-entropy method as in Section \ref{S6.3}, or to
adapt the ODE approach to non-positive but less singular
approximations (that we actually intend to do).


\subsection{Compacton for a cubic fifth-order NDE}

 For introducing a new model, unlike most of previous cases (excluding
 (\ref{HD1N}) in Section \ref{S9.3}), without any hesitation, we consider the cubic
 NDE--5
 \beq
 \label{cu1}
u_t=-(u^2 u_x)_{xxxx}+ u^2 u_x \quad \mbox{in} \quad \re \times
\re_+.
 \eeq
 We take the following
TW compacton with the specially chosen  wave speed $\l=-\frac 13$:
 \beq
 \label{cu2}
  \mbox{$
  u_{\rm c}(x,t)= f(y), \quad y=x+\frac 13 \, t
  \quad \Longrightarrow \quad -(f^3)^{(4)}+f^3=f \,\,\,\, \mbox{in}\,\,\,\, \re.
 $}
  \eeq
We next perform
 the natural change leading to a simpler semilinear ODE,
  \beq
  \label{cu3}
  F=f^3\quad \Longrightarrow \quad F^{(4)}=F-F^{\frac 13} \,\,\,\, \mbox{in}\,\,\,\, \re.
 \eeq
 This easy looking equation admits a nontrivial countable set of various
 compactly supported solutions that are analyzed by variational
 methods based on Lusternik--Schirel'man category and Pohozaev's fibering
 theory, \cite{GMPSob}. Here we stress our attention to
 the primary facts that are connected with the proposed concepts
 of entropy solutions.

\begin{figure}
\centering
\includegraphics[scale=0.65]{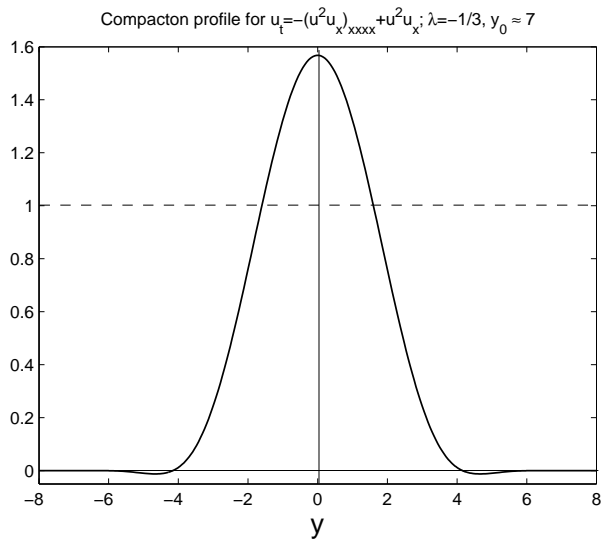}
 \vskip -.3cm
\caption{\small The first compacton profile of the ODE
(\ref{cu2}).}
\label{FC1}
\end{figure}

 The first and simplest compacton solution of the ODE (\ref{cu2})
 is shown in Figure \ref{FC1} that was obtained numerically with
 the  tolerances and regularization parameters
  $$
  {\rm Tols} = 10^{-10}
   \quad \mbox{and}
   \quad F^{\frac 13} \mapsto (\nu^2+F^2)^{-\frac 13}F \,\,\,
   \mbox{with also}
   \,\,\, \nu=10^{-10}.
 $$

\subsection{Oscillatory structure near interfaces: periodic
orbits}

\begin{figure}
\centering
\subfigure[zero structure $\sim 10^{-3}$]{
\includegraphics[scale=0.52]{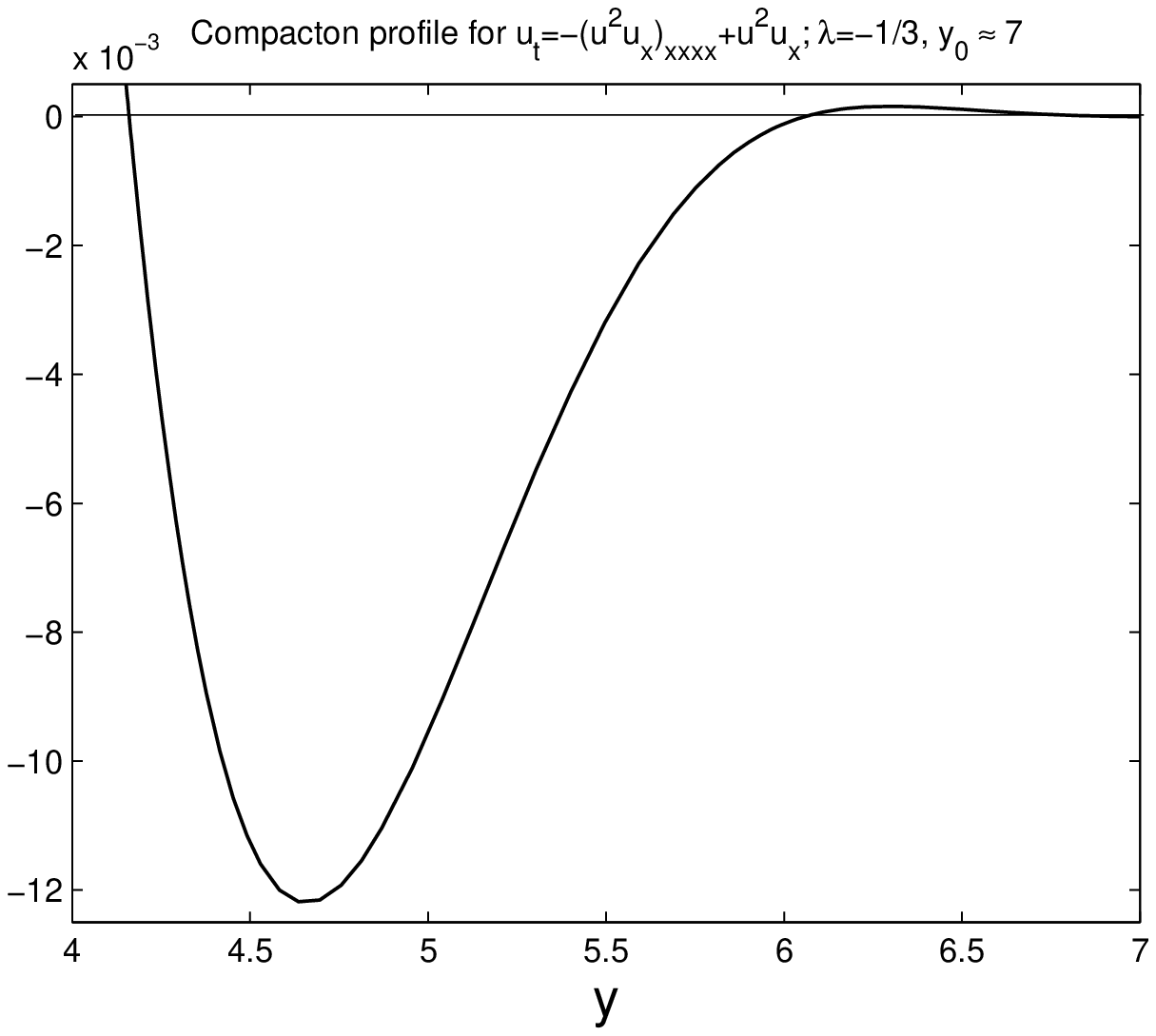} 
}
\subfigure[zero structure $\sim 10^{-5}$]{
\includegraphics[scale=0.52]{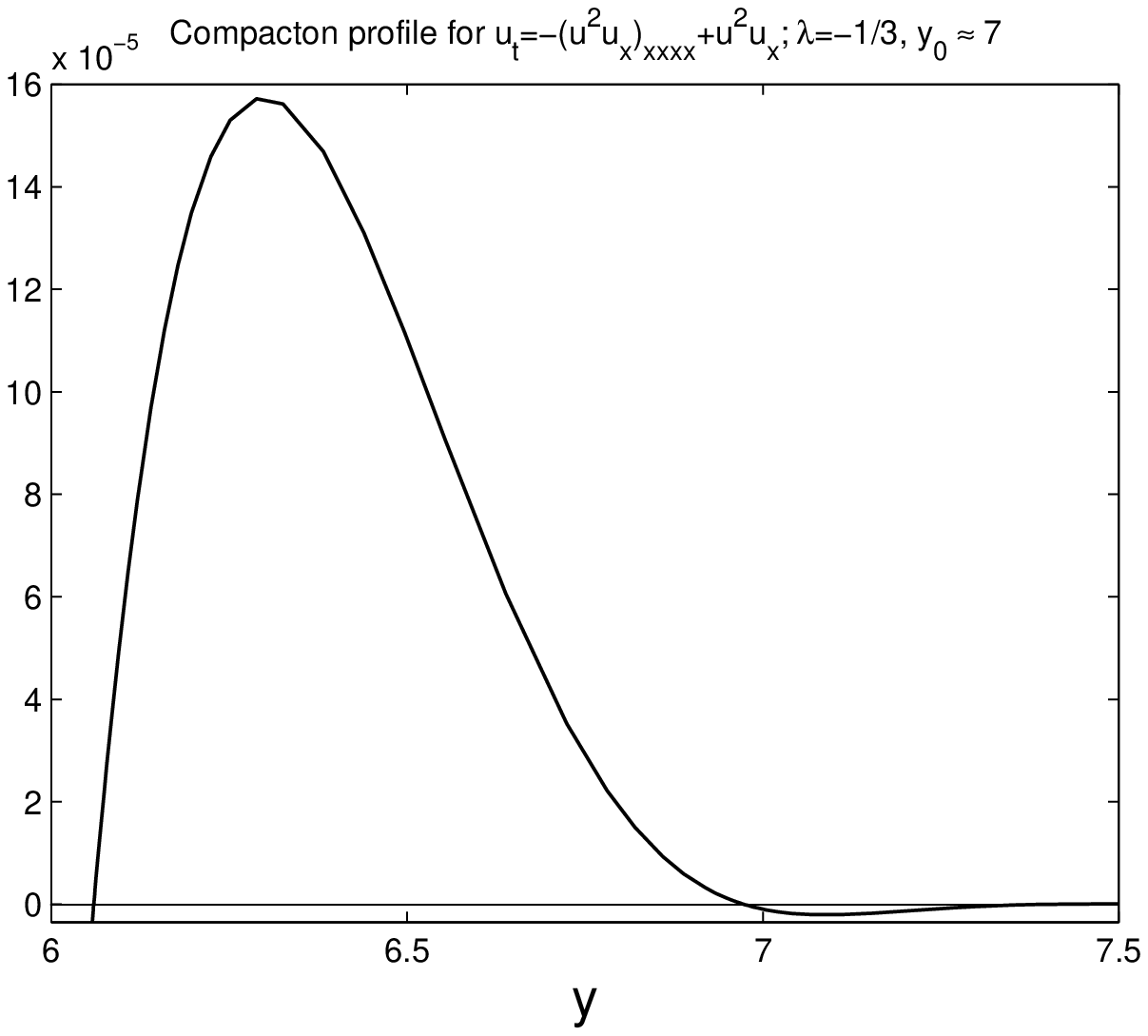} 
}
 \vskip -.4cm
\caption{\rm\small Enlarged zero structure from Figure \ref{FC1}
near the right-hand interface.}
 \label{FSign1}
\end{figure}

In general, it looks that this compacton profile does not differ
from those considered before as the explicit solutions in
(\ref{62}) or (\ref{789}). However,  there is a fundamental
difference that changes the mathematics of such solutions: for the
fifth-order NDE (\ref{cu1}), the profiles $f(y)=F^{\frac13}(y)$
are oscillatory and are of changing sign near finite interfaces.
 In Figure \ref{FSign1}, we show first three zeros near the
 interface at $y=y_0>0$ of the compacton profile from Figure
 \ref{FC1}.


In order to describe key features of such oscillatory behaviour at
the right-hand interface, as $y \to y_0>0$, when $F(y) \to 0$, we
perform an extra scaling by setting in the two leading terms of
the ODE
 \beq
 \label{cu4}
 F^{(4)}=-F^{\frac 13} \quad \Longrightarrow
 \quad F(y)=(y_0-y)^6 \varphi(s), \quad s=  \ln(y_0-y),
  \eeq
  where the {\em oscillatory component} $\varphi(s)$
  solves the following ODE:
   \beq
   \label{cu5}
   \begin{matrix}
   P_4(\varphi)
   \equiv {\mathrm e}^{-2s}
   \big[{\mathrm e}^{-s}({\mathrm e}^{-s}({\mathrm e}^{-s}({\mathrm e}^{-s}
   ({\mathrm e}^{6s}\varphi)')')')'\big]
    \ssk\ssk\\
   \equiv\, \varphi^{(4)} + 16 \varphi''' + 119\varphi''
 + 342 \varphi'
 + 360 \varphi=-\varphi^{\frac 13}.
 \end{matrix}
  \eeq
It turns out that the oscillatory behaviour near the interface at
$y=y_0^-$ (i.e., at $s= -\infty$) is given by a periodic solution
$\varphi_*(s)$ of the ODE (\ref{cu5}). Namely, we list the
following properties that lead to existence of a periodic orbit of
changing sign:

 \begin{proposition}
 \label{Pr.Per2}
 The fourth-order dynamical system $(\ref{cu5})$ satisfies:

\noi {\rm (i)} no orbits are attracted to infinity as  $s \to +
 \infty$;

\noi {\rm (ii)} it is a dissipative system with a bounded
absorbing
 set; and


\noi {\rm (iii)} a nontrivial periodic orbit $\varphi_*(s)$
exists.

\end{proposition}

\noi{\em Proof.} (i) The operator in (\ref{cu5}) is asymptotically
linear \cite[p.~77]{KrasZ} with the derivative at the point  at
infinity $P_4$ that has the characteristic equation
 $$
 p_4(\l)=\l^4 + 16 \l^3 + 119\l^2
 + 342 \l
 + 360 \equiv (\l+6)(\l+5)(\l+4)(\l+3)=0.
  $$
Therefore, all eigenvalues are real negative, $-6$, $-5$, $-4$,
and $-3$, so infinity cannot attract orbits as $s \to + \iy$. (ii)
is a corollary of (i) after an extra scaling.


(iii) Existence of a periodic orbit for dissipative systems is a
standard result of degree theory; see \cite[p.~235]{KrasZ}. We
complete the proof of existence by using a shooting argument as in
\cite[\S~7.1]{Gl4}.  $\qed$


It turns out that the periodic solution $\varphi_*(s)$ is
exponentially stable as $s \to +\infty$ (this is not easy to see
from the ODE (\ref{cu5}) by linearization and interpolation of the
third term with the coefficient 119). The asymptotic stability of
this periodic orbits is illustrated in figures in
\cite[p.~187]{GSVR}.

Thus, at the singular end point $y=y_0^-$, the dynamical system
(\ref{cu3}) generates a two-dimensional bundle of orbits with the
behaviour
 \beq
  \label{as11}
 F(y) = (y-y_0)^6[\varphi_*(s+s_0)+...], \quad y_0 \in \re_+,
 \,\,\, s_0 \in \re,
 \eeq
 where $s_0$ is an arbitrary  phase shift of the periodic motion.
 Thus, the interface point $\{y=y_0^-, \, F=F'=F''=F'''=0\}$ is a  complicated
 singular point (a zero) of the dynamical system (\ref{cu3}), so
 one needs to check whether it corresponds to an entropy
 solution-compacton.
  It is worth mentioning that the
  2D bundle (\ref{as11})  matches with precisely two symmetry conditions at
 the origin,
  $
  F'(0)=F'''(0)=0,
  $
and  the existence of the compacton is confirmed by  variational
methods, \cite[\S~5]{GMPSob}.

\subsection{Compactons are $\d$-entropy: a formal illustration}

The oscillatory  behaviour (\ref{as11}) of the
compacton near  finite interfaces makes impossible to use the
{\em positive} analytic $\d$-approximation as for the NDE--3  in
(\ref{63}). Indeed, the same procedure for $F=f^3$ now leads to
the ``regularized" ODE
 \beq
 \label{rr1SS}
 F_\d:\quad
 F^{(4)}=F-F^{\frac 13}+ C_\d, \quad \mbox{where}
 \quad C_\d=\d^{\frac 13}-\d, \quad \mbox{so}
 \quad F_\d(y)\to \d>0, \,\,\, y \to \infty.
  \eeq
This gives the family $\{F_\d\}$ consisting of  functions
$F_\d(y)$ that change sign {\em finitely many times} for all
sufficiently small $\d>0$. These  approximations $F_\d$ are less
singular than the limit compacton profile $F(y)$, which according
to (\ref{as11}) is infinitely oscillatory as $y \to y_0^-$.

\begin{figure}
\centering
\subfigure[$F_\d \to F$]{
\includegraphics[scale=0.52]{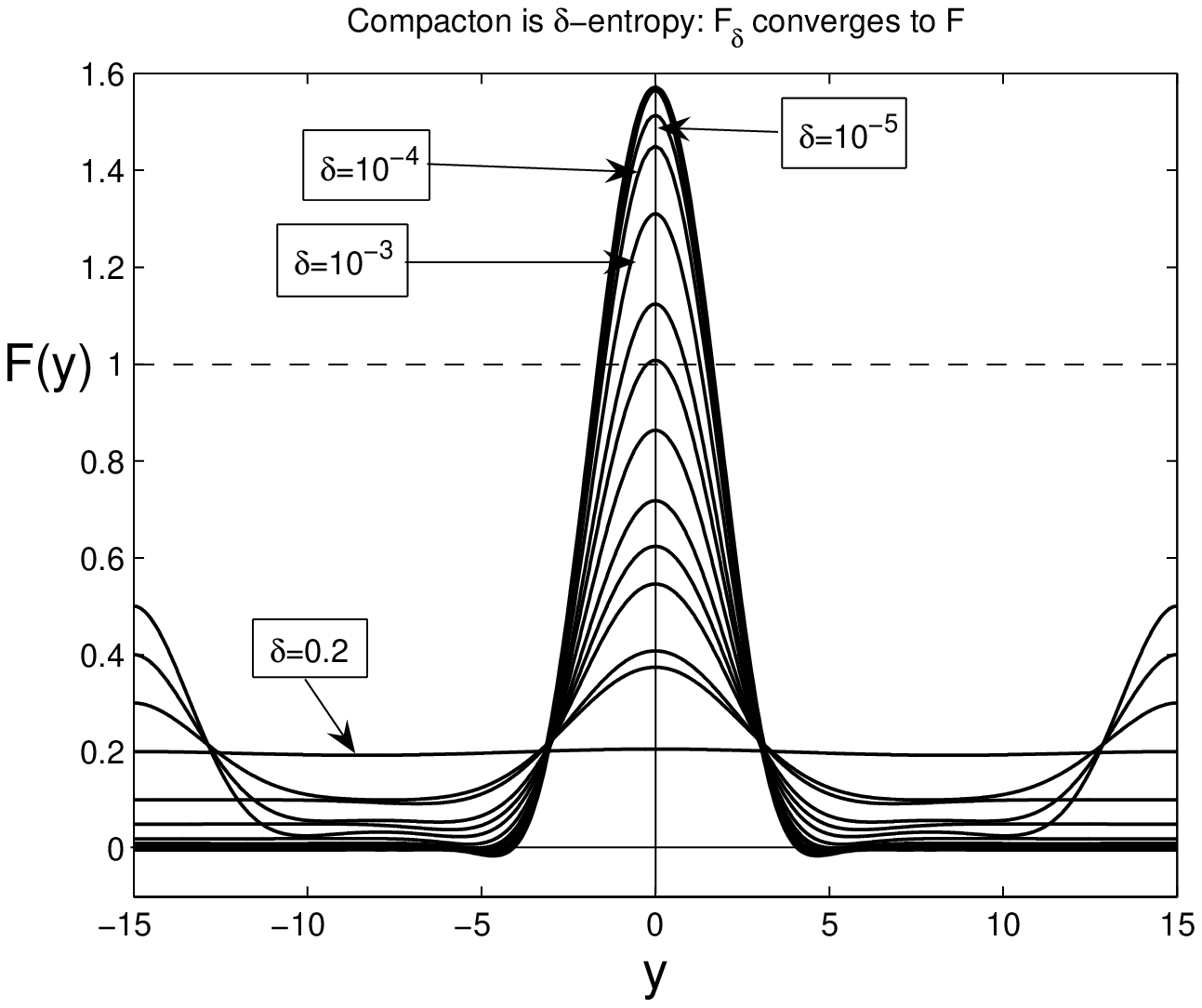} 
}
\subfigure[zeros of $F_\d$, enlarged]{
\includegraphics[scale=0.52]{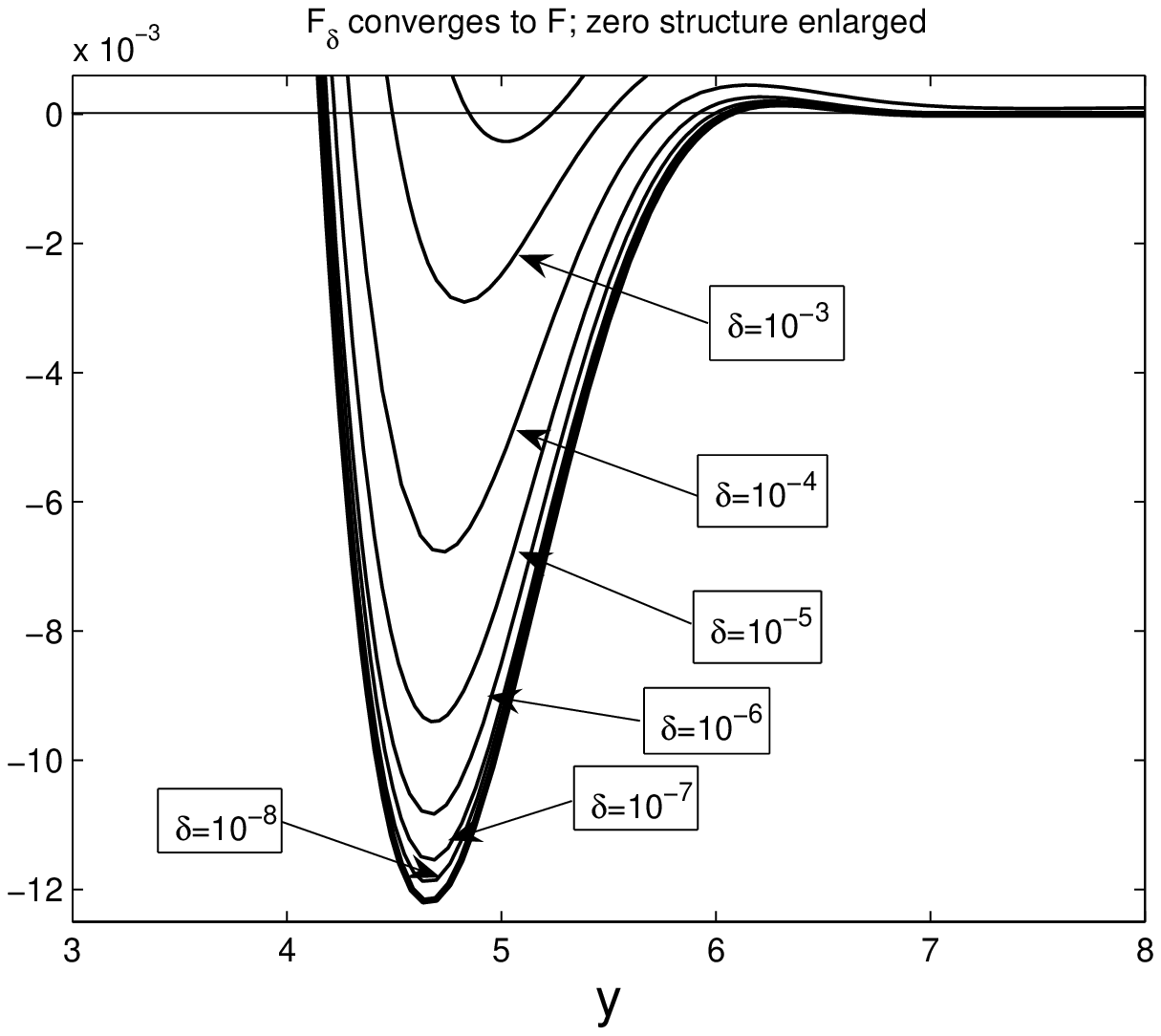} 
}
 \vskip -.3cm
 \caption{\small Convergence of the regularizing sequence
$\{F_\d\}$ from (\ref{rr1SS}) to the compacton $F(y)$ (a) and
enlarged finite zero structure of $F_\d$ for small $\d>0$ (b).}
 \label{FD1}
\end{figure}

The solvability of the approximating problem (\ref{rr1SS})
can be traced out by the same variational method.
 Then the convergence
\beq
 \label{Conv1}
  F_\d(y) \to F(y) \quad \mbox{as}
  \quad \d \to 0^+ \quad \mbox{uniformly}
  \eeq
 is associated with the
 stability of critical values of functionals; see
 \cite[p.~387]{KrasZ}.
 This
$\d$-approximation  is shown in Figure \ref{FD1}(a), where the
convergence (\ref{Conv1})
 is rather slow and is observed starting from $\d=10^{-3}$ only, with the accuracy
 about $0.2$.
 For $\d=10^{-2}$, the approximating profile $F_\d(y)$ is still
 almost four times less than $F(y)$ at the origin. The accuracy
 $0.1$ is achieved for $\d=10^{-5}$. In (b), up to $\d=10^{-8}$, we show the zero
 structure of $F_\d(y)$ close to $y_0$, which, since $F(y) \approx
 \d>0$ for $y \gg 1$, is  finite and each zero is
 transversal. These confirm that the approximating sequence
 $\{F_\d\}$, though is of changing sign, is less singular than the
 compacton profile $F(y)$ itself.


An alternative approximating approach of such compactons
 is developed in \cite{GMPSobII}, where $F$ is approximated as $\e \to 0^+$ by
 the analytic family $\{F_\e\}$ of solutions of the regularized ODE
  $$
  F_\e: \quad F^{(4)}=F-(\e^2+F^2)^{-\frac 13}F \quad (\e>0).
  $$
Incidentally, this approach makes it possible to trace out the
{\em Sturmian index} of some solutions by a homotopic connection
to variational problems with known ordered set of critical points
and known number of zeros for each of them, \cite{GMPSob}.


\enddocument